\documentclass[12pt]{paper}
 \usepackage{harvard}

  \usepackage{cite}
 \usepackage{amsmath}  
\usepackage{amsfonts}
\usepackage{amssymb}
\newcommand{\Smal}{{\bf (Small)}\ \ }
\newcommand{\Comp}{{\bf (Comp)}\ }
 
 \usepackage{cleveref}
 \usepackage{cleveref}
 \usepackage{centernot}
 \usepackage{epsfig}
\usepackage{graphics}
 \usepackage{color}

\usepackage{color}
\usepackage{epsfig,cmmib57}
\usepackage{amssymb,amsmath,exscale}
 \numberwithin{equation}{section}

\newcommand{\ZEP}{\epsilon}

\newcommand{\ZOMq}{\Omega}

\newcommand{\zg}{\gamma}

\newcommand{\intt}{\int_0^t}
\newcommand{\ints}{\int_0^s}

\newcommand{\GP}{{\bf (GP)}\ }
\newcommand{\MC}{{\bf (MC)}\ }
\newcommand{\J}{{\bf (J)}\ }
\newcommand{\CG}{{\bf (CG)}\ } 

\newcommand{\ZOM}{\omega}
\newcommand{\zaa}{\alpha}

\newcommand{\zt}{\tau}

\newcommand{\ZLA}{\label}
\newcommand{\ZIN}{\infty}

\newcommand{\zzr}{{\rm I\hskip-2.1pt R}}

\newcommand{\ZD}{\;\mbox{\rm d}}
\newcommand{\zl}{\lambda}
\newcommand{\ZSI}{\sigma}

\author{
L. Pandolfi\thanks{Retired from the Dipartimento di Scienze Matematiche ``Giuseppe Luigi Lagrange'', Politecnico di Torino, Corso Duca degli Abruzzi 24, 10129 Torino, Italy (luciano.pandolfi@polito.it)}
}
\title{Linear systems with persistent memory: An overview of the biblography on controllability\thanks{
This papers fits into the research program of the GNAMPA-INDAM and has been written in the framework of the   ``Groupement de Recherche en Contr\^ole des EDP entre la France et l'Italie (CONEDP-CNRS)''.}}
\begin{document}
 
 \maketitle 
 
 \begin{abstract}
 We discuss the bibliography on control problems for systems with persistent memory
 \end{abstract}

\section{Introduction}
 
We present an overview of the  controllability results for
 a class of systems described by equations whose prototype has the following form:
 \begin{equation}\ZLA{eq:memo}
w'=\alpha A  w(s)+\intt N(t-s) A w(s)\ZD s  +F(t)\,.
 \end{equation}
 Here $A$ is the infinitesimal generator of  a $C_0$-semigroup in a Hilbert space $X$.

 Actually, in the cases   that have been studied from the point of view of controllability $A$ is the laplacian or the Navier operator of three-dimensional elasticity (possibly a second order uniformly elliptic operator) when studying wave type equation or the bilaplacian (in the study of the vibration of a plate).
 
 Of course we must associate suitable initial and boundary conditions, as described below.

 The function $N(t)$ is called the \emph{memory } or the \emph{relaxation} kernel. We assume
 $N\in   C^k([0,+\ZIN))$ ($k\geq 1$ as specified below); $\zaa\geq 0$ and  if $\zaa=0$ then we assume $N(0)>0$ (these   conditions imply well posedness of the process, see below).

When $\zaa=0$ and $N(0)>0$, a case encountered in particular   in the study of viscoelasticity, \cref{eq:memo} is often written in the following form form:
\begin{equation}
\ZLA{eq:memoSecondOrder}
w''= N(0) Aw+\intt M(t-s)Aw(s)\ZD s +G(t)\,,\qquad M=N'\,.
\end{equation}

There is a large recent literature on controllability of fractional systems (in this case $N(t)$ is weakly singular at $t=0$) but we are not going to consider this case. We refer for example to~\cite{BonaccorsiJEVOLEQ2012,WangFeskanEVEQCONTRTH2017,JacobBirgitJFA2010,SaktivelCOMPUTERS2006,SaktivelJMAA2007}
 and reference therein.  
 We note also that we confine ourselves to review controllability of linear systems.    
 
  \section{A short overview of the origin of the models}

It might have an interest to shortly describe the origin of this class of systems which have been introduced in the second half of the $XIX$ century independently by several authors and for different reasons.   One of the earliest appearence  of systems of the form~(\ref{eq:memo}) is in a paper of J.C. Maxwell, 
which uses \cref{eq:memo} to model a viscoelastic fluid (see~\cite{MaxwellPHILOSOPHICAL1867}). In that paper J.C. Maxwell introduced a device which then became   standard in the engineering description of viscoelastic materials: any element of a viscoelastic fluid is described as a ``cell'' which contains both a spring and a damper in series, so that the stress   solves
 \[
\ZSI'(x,t)= -\frac{1}{\zt}\ZSI(x,t)-\kappa \nabla w(x,t) \qquad \zt>0\,,\ \eta>0 \,.
 \]
 Combining this material law with conservation of momentum, we get \cref{eq:memo} with $ \zaa=0 $ and a memory kernel of the form $ N(t)=\beta e^{-\eta t} $. If instead the element is modelled as a spring and damper in parallel then we get Voigt   model (the case $ \zaa>0 $ is a generalization of Voigt  model). A nice discussion is in~\cite[Ch.~5]{KolskyBOOK} and in~\cite[Ch.~1]{JosephLIBRO1990}.

The use of ``general'' kernels in viscoelasticity appears in the linearization of the  Boltzmann equations and are justified by the \emph{Boltzmann superposition principle} (see~\cite{BoltzmannWIN1874,BoltzmannWIED1878}), which is a consequence of the following observation:   if a viscoelastic body is subject to a certain stress, it will need some time to attain the corresponding deformation: in a certain small interval of time $\Delta s$ it will only undergone a ``partial'' deformation. Conversely, when the external stress is removed, the deformed configuration and the corresponding interior stress will slowly relax. The strain in position $x$ approximately equal to  $  \nabla w (x,s)$ on an interval of time $(s,s+h)$ produces a stress at the future times and at the time $t>s  $  it contributes    a stress approximately equal to  $-N(t-s)\nabla w (x,s) h$, at the same point $x$. At time $t$ the stresses produced by the previous strains $\nabla w(x,s)$    are summed, so leading to an equation of the form~(\ref{eq:memo}).  These ideas, further developed by 
V. Volterra, see~\cite{VOLTERRArendiconti1909,VolterraACTAMATHEM1912,VolterraLIBROfUNZINEA}, are subsumed by the following law which relates strain and stress, and which replaced Hook law (for linear systems):
\[ 
\ZSI(x,t)=-\nabla\left (\zaa w(x,t)+\intt N'(t-s) w(x,s)\ZD s\right )\,.
 \]
 
%
%
As first proved in~\cite{KOHLRAUSCH1863}, 
 are viscoelastic materials.  At the time, glasses where used as dielectric in capacitors (Leiden Jars). On the suggestion of J.C. Maxwell, 
the viscoelastic properties of glasses where exploited in~\cite{HopkinsonLONDON1877,HopkinsonLONDON1876} to
explain the relaxation of the potential of a capacitor when the dielectric is not the void. In these papers, J. Hopkinson realized that  kernels $N(t)=\sum \beta_k e^{-\eta_k t} $ satisfactorily fit experimental data and this fact was interpreted as follows: 
a glass constituted by \emph{only one} silicate  would produce  a discharge law which relaxes in time with \emph{one exponential  } $ \beta  e^{-\eta  t}$. Real glasses
 are a mixture of different silicates. Each one of them produces a relaxed discharge,  relaxing with its own exponential kernel $\beta e^{-\eta t}$: these different exponential relaxation kernels have to be linearly combined in order to have an accurate description of the process. This justifies the appearance of the  kernel  
\[
N(t)=\sum \beta_k e^{-\eta_k t}
\]
which is called a  \emph{Prony sum} in engineering literature.

  Few years later, in 1889, M.J. Curie in his studies on  piezoelectricity   noted  that the discharge of capacitors with certain crystals as dielectric  is best described by a rational law $i(t)={\rm (const)}/t^\zg$ (see~\cite{CurieANNALES1889}). This leads to a fractional integral law for the current released by a capacitors with cristalline dielectric (see~\cite{WesterlundTRANSDIEL1994}). A fractional integral is a common model for distributed systems with persistent memory, and now much studied also from the point of view of controllability.   
As already stated however we are not going to review results concerning singular kernels.

H. Jeffreys introduced the special case $\zaa>0$ and $N(t)=\beta e^{-\eta t}$ in order to take into account viscosity of the ocean water (see~\cite{Jeffreys1924}).

 Finally we mention that     \cref{eq:memo} is also used
  in the study of   ferromagnetic materials, as first proposed in by D.~Graffi in~\cite{GraffiCIMENTO1928,GraffiLOMBARDO1936} (see also~\cite[Ch.~3]{VolterraACTAMATHEM1912}.

Note that when  $N(t)=\beta e^{-\eta t}$ and $\zaa=0$ then \cref{eq:memo} (with $A=\Delta$ and $F(t)=0$) is the integrated version of the telegrapher's equation
\begin{equation}\ZLA{eq:tele}
w''=-\beta w'+\eta\Delta w
\end{equation}
and it is known that solutions of \cref{eq:tele} have a wavefront, a discontinuity in the direction of the propagation, but not backward (unless $\eta=0$ when the equation reduces to the wave equation). This ``forward wavefront'' is observed in the diffusion of solutes in solvents with complex molecular structure, and for this reason ~\cref{eq:tele} and more in general (\ref{eq:memo}) are used also to model nonfickian 
diffusion (see~\cite{ChristensenBOOK1982}), for example in biology. It seems that the use of \cref{eq:memo}
in order to model diffusion in polymers goes back to the late forthies of the $XX$-century (see~\cite{DEkEELIUHinestroza}), but I was unable to locate the earliest appearence of this model.

\Cref{eq:tele} was first introduced in thermodynamics in~\cite{CattaneoMODENA1949} as an ``hyperbolic'' model for heat diffusion with a finite propagation speed.  Finally,
\cite{ColemanGurtinZEITSCH1967,GurtinPipkinARCHIVE1968}
  introduced the general   \cref{eq:memo} to model heat processes (in fact, the model introduced in~\cite{GurtinPipkinARCHIVE1968} is more complex, with two memory kernels, but it can be reduced to  \cref{eq:memo}   if the memory kernels are known).
 
 The material property which leads to \cref{eq:memo} in thermodynamics is the following temperature/flux law, which replaces Fourier law:
 \begin{equation}\ZLA{eq:Flusso}
 q(x,t)=-\zaa\nabla w(x,t)-\intt N(t-s)\nabla w(x,s)\ZD s 
  \end{equation}
  (so that in this case  $ A=\Delta $).

Stimulated by so many applications,  \emph{distributed system with persistent memory} described by \cref{eq:memo} (or more general equations) with suitable initial and boundary conditions, became a subject of intense study (see for example the
 books~\cite{Amendola2012,ChristensenBOOK1982,FabrizioMorroCBMS,RenardyHrusaLibro1987,RenardyCBMS2000}). Among the problem which received the wider attention we list:
\begin{itemize}
\item definition of energies;
\item conditions   imposed to the kernel by fundamental physical law (dissipation of energy and the second principle of thermodynamics);
\item conditions on the kernel which ensure   (strong or exponential) stability; 
\item conditions on the   kernel such that signal propagates with finite or infinite speed;
 
\item identification problems (of inputs or coefficients, in particular of the relaxation kernel) on the basis of suitable measurements taken on the system;
\item control problems, which is the most recent topic.
 
 \end{itemize}

Papers on the finite or infinite speed of propagation are particularly important to us, since they relate \cref{eq:memo} either to the heat or the wave equation. Among the relevant paper we cite~\cite{DeschINTEGRALEQ1985,FisherGurtin1965,HerraraGurtin1965,NarainJosephLINEARIZED1983} and in particular~\cite{RenardyREOLOGICAL1982}.
It is proved in this paper that when $\zaa=0$ signals propagate with a finite speed if $N(t)$ is smooth for $t\geq 0$ and if $N(0)>0$.    (the velocity of propagation is $\sqrt{N(0)}$). In this case the initial conditions are not smoothed.   If instead   $N(t)\sim 1/t^\zg$, $\zg\in(0,1)$, like in Curie law, then the velocity of propagation is still finite but the signal is smoothed. In particular this implies that  exact reachability of an arbitrary specified target   is impossible in this case.

A rich overview on \cref{eq:memo}  is in~\cite{JosephPreziosi1989,JosephPreziosiAddendum1990}.

  \subsection*{Classification of systems}
  
  As we stated, we are not going  to consider the case that $N(t)$ has a singularity at $t=0$ and the previous observations suggest that we distinguish between the following two classes of systems:
  \begin{itemize}
\item[\CG] (after Colemann and Gurtin) is the case $\zaa>0$. The special case $N=0$ is the heat equation. If $\zaa>0$ and $N(t)=1$ we get the Voigt  model.
\item[\GP] (after Gurtin and Pipkin) is the case $\zaa=0$ and $N(0)>0$. The special case $\zaa=0$ and $N(t)\equiv 1$ is the wave equation.
\end{itemize}
  
 Special cases are
  \begin{itemize}
  \item[\MC] (after Maxwell and Cattaneo) which is \GP\ when   $ N(t)=\beta e^{-\eta t} $.
  \item[\J] (after Jeffrey) which is \CG with  $ N(t)=\beta e^{-\eta t} $.
  \end{itemize}

There is a large literature on controllability of the telegrapher's or Voigth equations. We are not going to review these results unless strictly relevant to \CG or \GP models (in the form \MC or \J\!\!).

  We note that the special case of \MC with $ F(t)=N(t)\xi =e^{-\eta t} \xi$ is the MGT equation of acoustics ($ \xi $ depends on the initial conditions), see~\cite{DelloroAMO2017,BucciPandolfiArxive}.

Finally:
\begin{itemize}
\item     if $\alpha=0$ and $N(0)=0$,  then 
\cref{eq:memo} might not be well posed, as seen from the following example:
\[
w'=\intt (t-s)\Delta(s)\ZD s
\]
which is the integrated version of
\[
w'''=\Delta w\,.
\]
This equation, for example on $(0,\pi) $ with Dirichlet homogeneous boundary condition, is not well posed, see~\cite[p.~99]{FattoriniBOOKcauchy1983}.
  \item
  The differential operators (denoted $ A $) inside and outside the memory integral might  be different, as in~\cite{VlasovJOURNMATHSCI}. 
 
  \end{itemize}

\section{Euristics on the relations of the wave and heat equations with the equation with persistent memory}

Finite propagation speed and the presence of the wavefront suggest  that certain properties of the wave equation are retained by  the \GP model.  In order to see the relations, let us look at the following figures.

We consider 
\[
w'=\zaa w_{xx}+\int_0^t N(t-s) w_{xx}(s) d s\,,\qquad  x\in (0,\pi) 
\]
and
\[
w(0,t)=w(\pi,t)=0\,,\qquad 
w(x,0)=\left\{
 \begin{array}{lll}
1 &{\rm if}&x\in ((1-\beta)\pi),(1+\beta)\pi)\\
0&& {\rm otherwise.}
 \end{array} 
\right.
\]
We choose
\[
N(t)=3e^{-t}-3e^{-2t}+e^{-4t}\,,\quad 
\beta=1/20\qquad \mbox{(note $N(0)=1$)}\,.
\]

\begin{figure}[h]
\caption{\label{Euris:1}The solution at the time $T=\pi/4$ of \GP compared with the solution of the wave equation  (left) and of \CG compared with the solution of the heat equation (right). } 
\vspace{0cm}
\begin{center}
 \includegraphics[width=5cm]{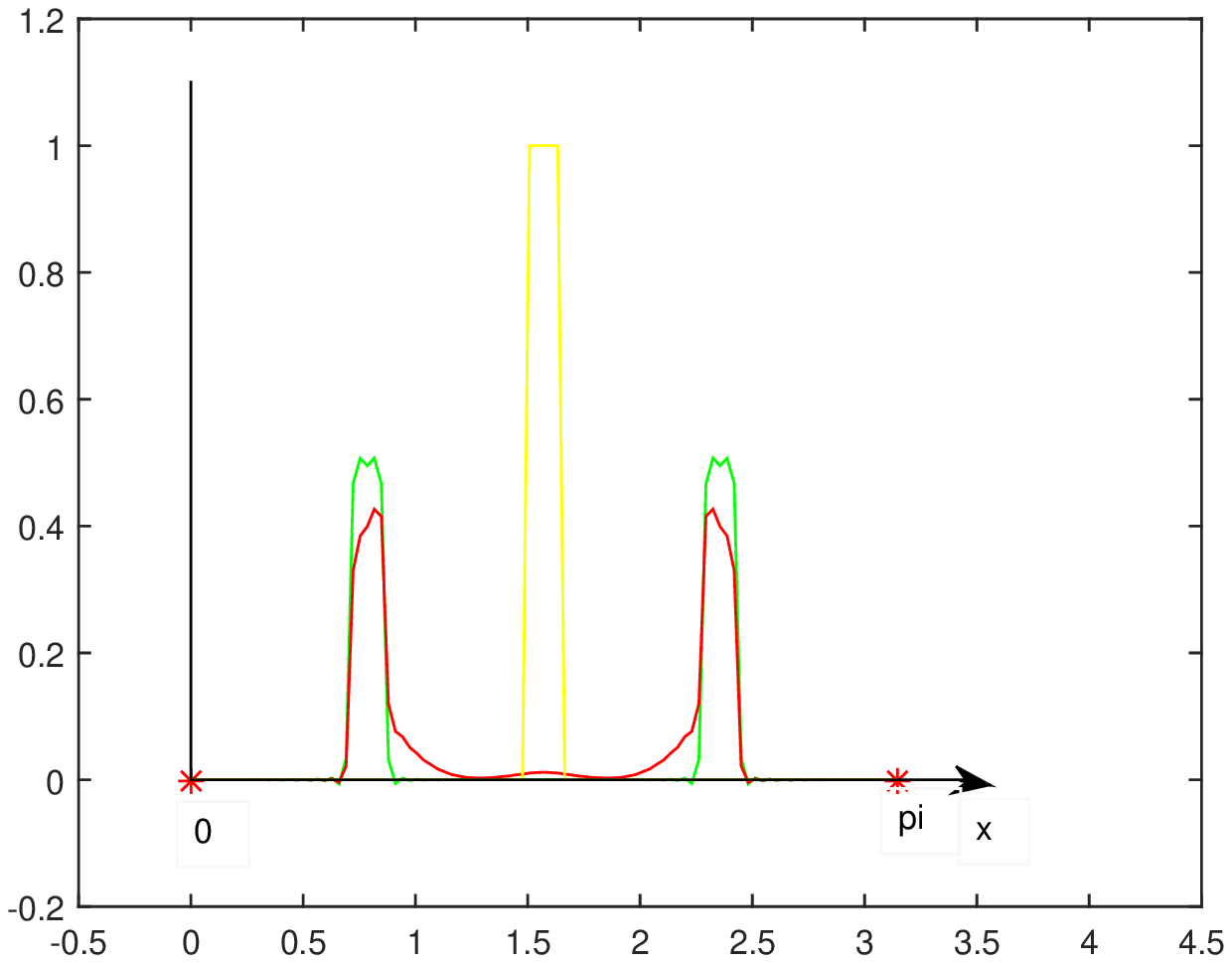} 
  \includegraphics[width=5cm]{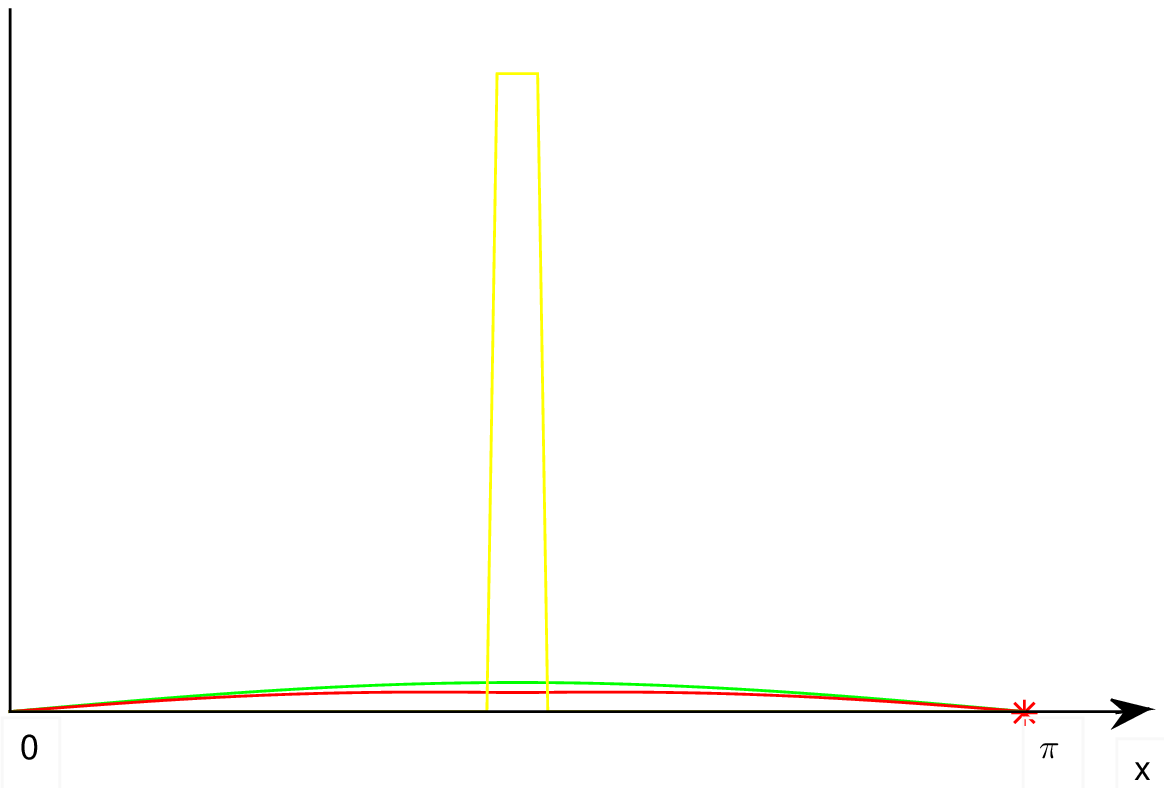}
\end{center}
\smallskip
\center{}
\end{figure}

\Cref{Euris:1} shows the plots at time $T=\pi/4$ of the solution of \GP (left) and of \CG (right, here $\zaa=1)$

\begin{center}
\begin{figure}[h]
\caption{\label{Euris:2}The solutions of the wave and \GP equations   after the reflection. Note that the solution of \GP is not zero on any interval.} 
\vspace{0cm}
\begin{center}
\includegraphics[width=7cm]{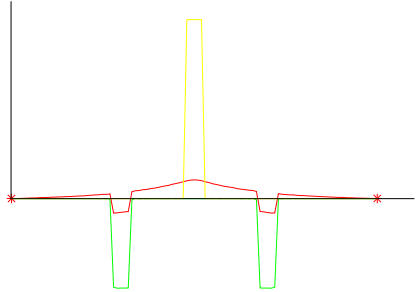}  
  \end{center}
\smallskip
\center{}
\end{figure}
\end{center}

The noticeable fact is that \CG is quite similar to the heat equation while \GP displays a significant similarity and also a significant difference with the wave equation. The velocity of the wave is $\sqrt{N(0)}=1$ so that at time $t=\pi/4$ the wave did not hit the ends of the interval and it displays the same wavefront in the direction of the motion as the wave equation but it does not became zero backward: the solution of the wave equation leaves no memory of itself after it has passed over a point $P$ while the memory of \GP persists. This is confirmed by \cref{Euris:2} which shows the
solutions \GP and of the wave equation after the reflection, i.e. at time $T>\pi/2$.

\begin{figure}[h]
\caption{\label{Euris:3}Low and high frequency content of the solutions at the time $T=\pi/4$ of the wave and \GP equation. } 
\vspace{0cm}
\begin{center}
\includegraphics[width=5cm]{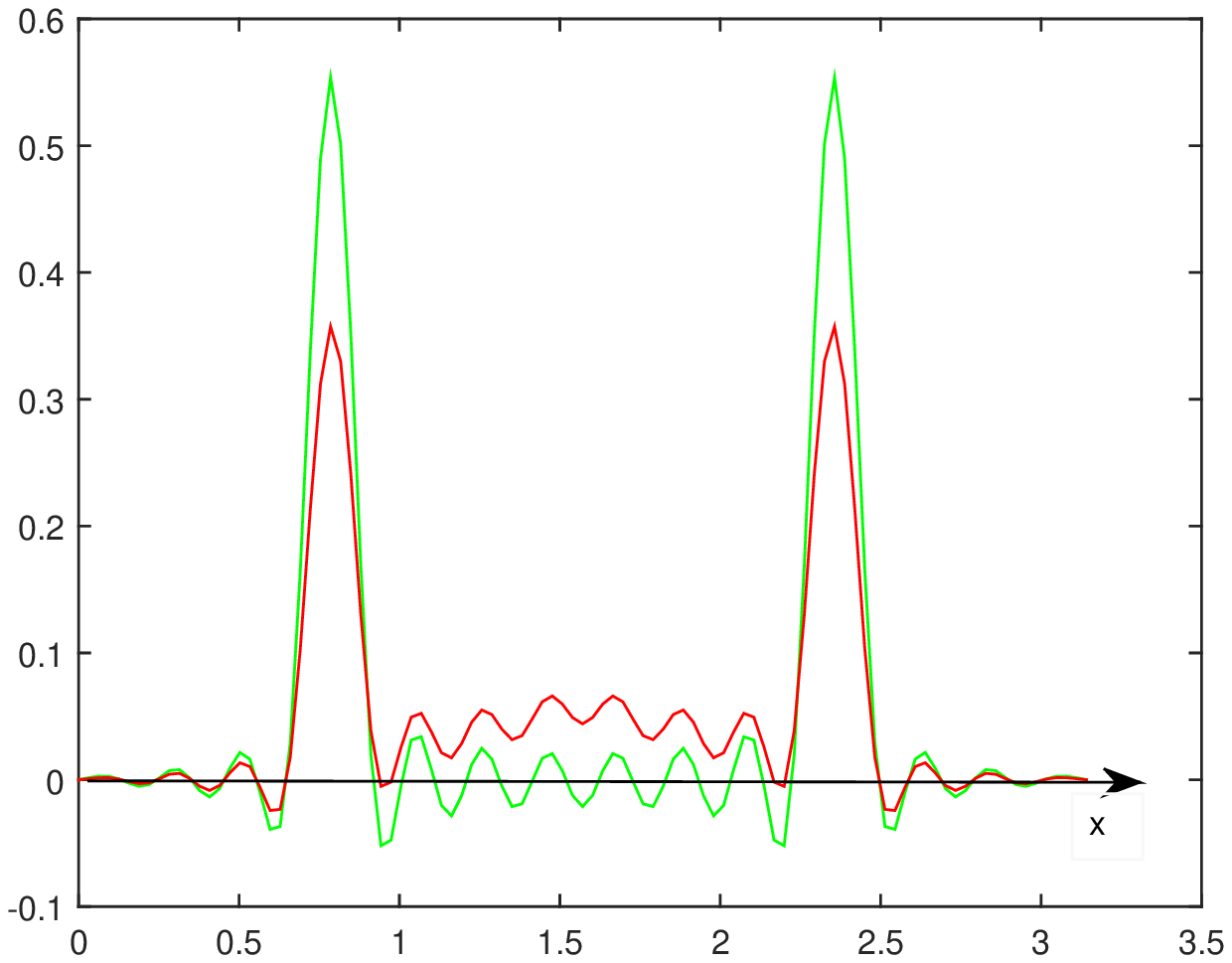}  
\includegraphics[width=5cm]{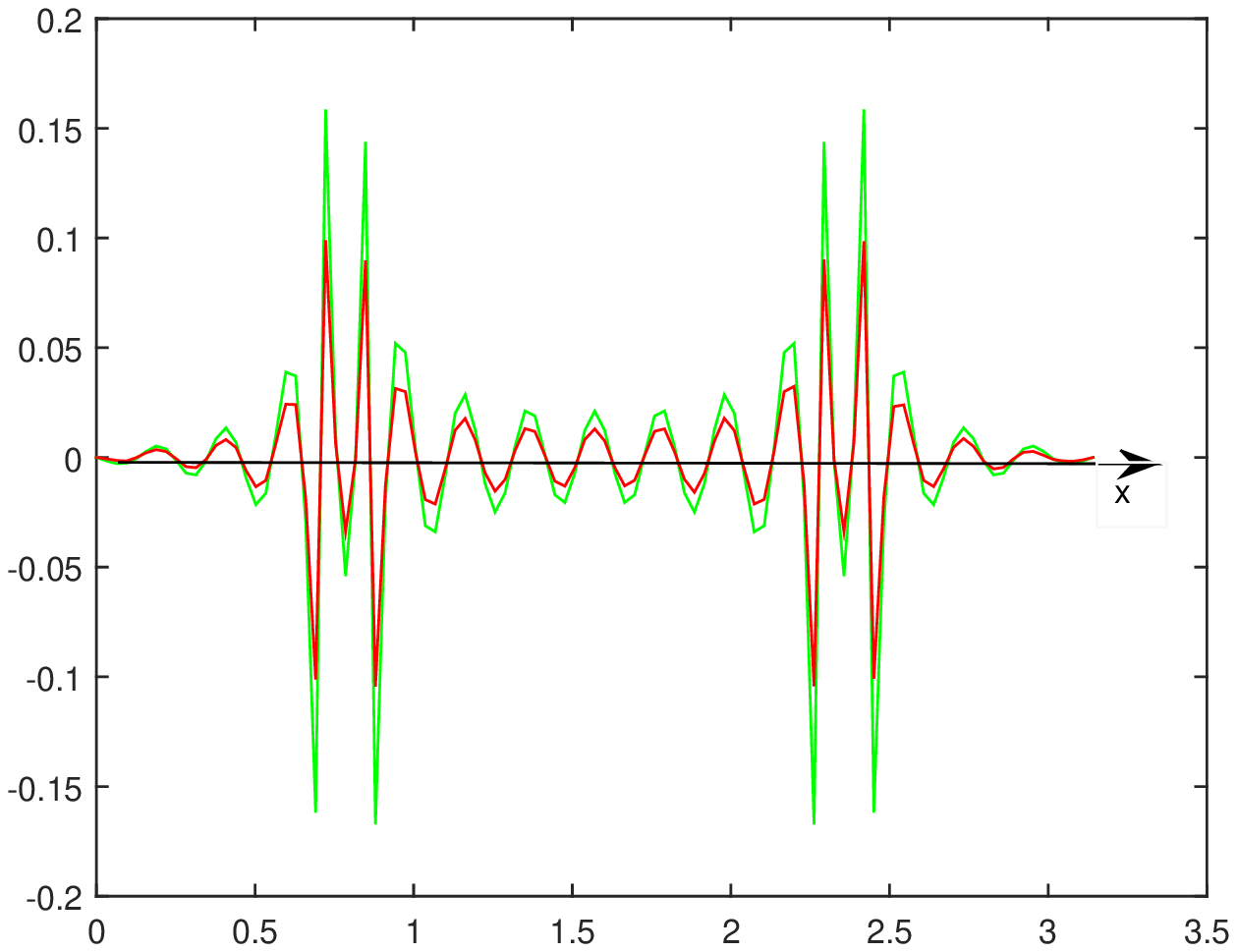}  
\end{center}
\smallskip
\center{}
\end{figure}

 It is clear that the most interesting features are displayed by \GP and it makes sense to try to understand better the discrepancy between the solutions of \GP and the wave equation. The next plots, in \cref{Euris:3}, decompose the solutions at the time $T=\pi/4$ (as displayed in \cref{Euris:1}) of \GP and of the wave equations in the contribution of the low and high frequencies (the frequency cut is   at the 30-th eigenvalue)

These figures display the well known fact that at high frequency a viscoelastic material behaves like a purely elastic one, and the discrepancy is mostly due to the    low frequency components.

%
%

\begin{center}
\begin{figure}[h]
\caption{\label{Euris:4}The solutions at time $T=3$ of \MC compared with the solutions of the wave equation (left) and of the heat equation (right).} 
\vspace{0cm}
\begin{center}
\includegraphics[width=4cm]{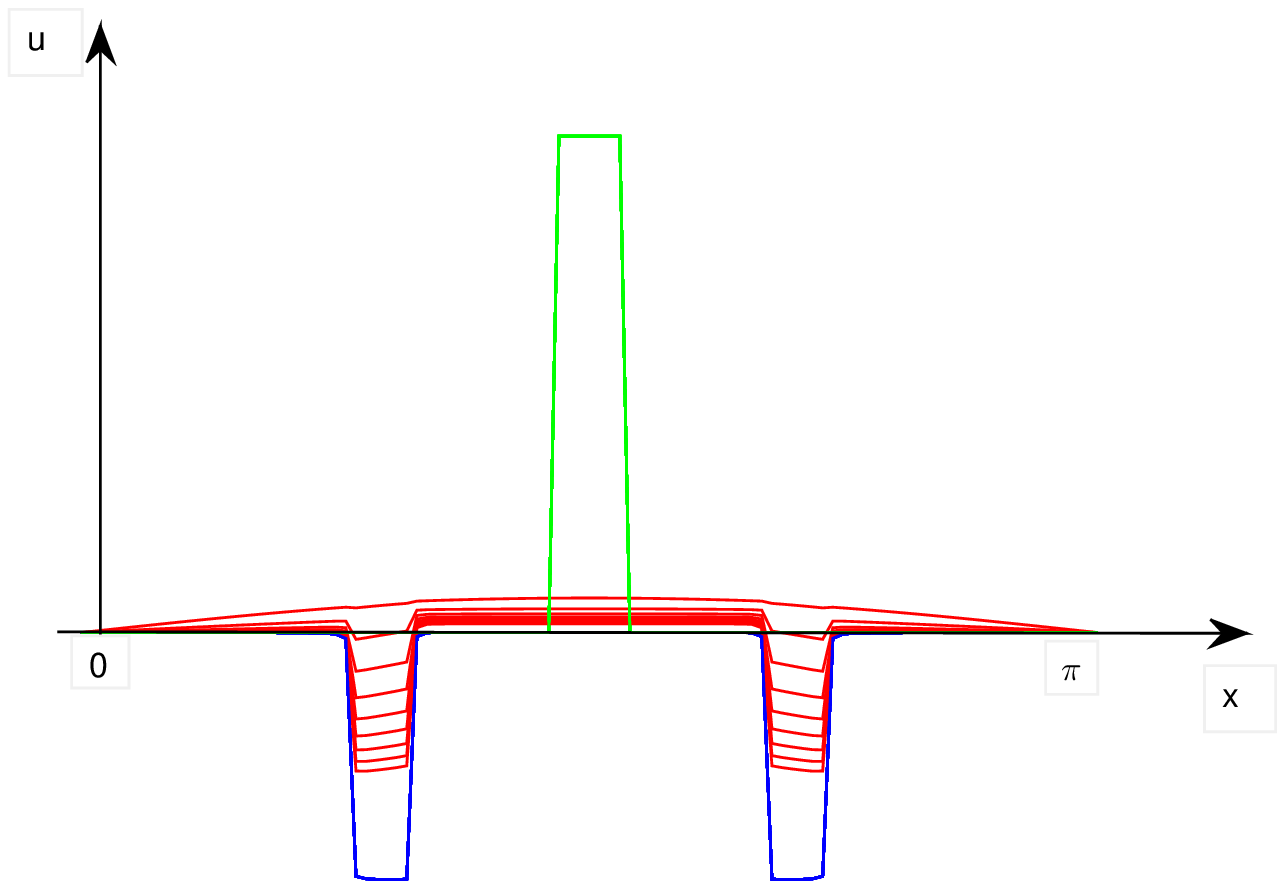} 
\includegraphics[width=4cm]{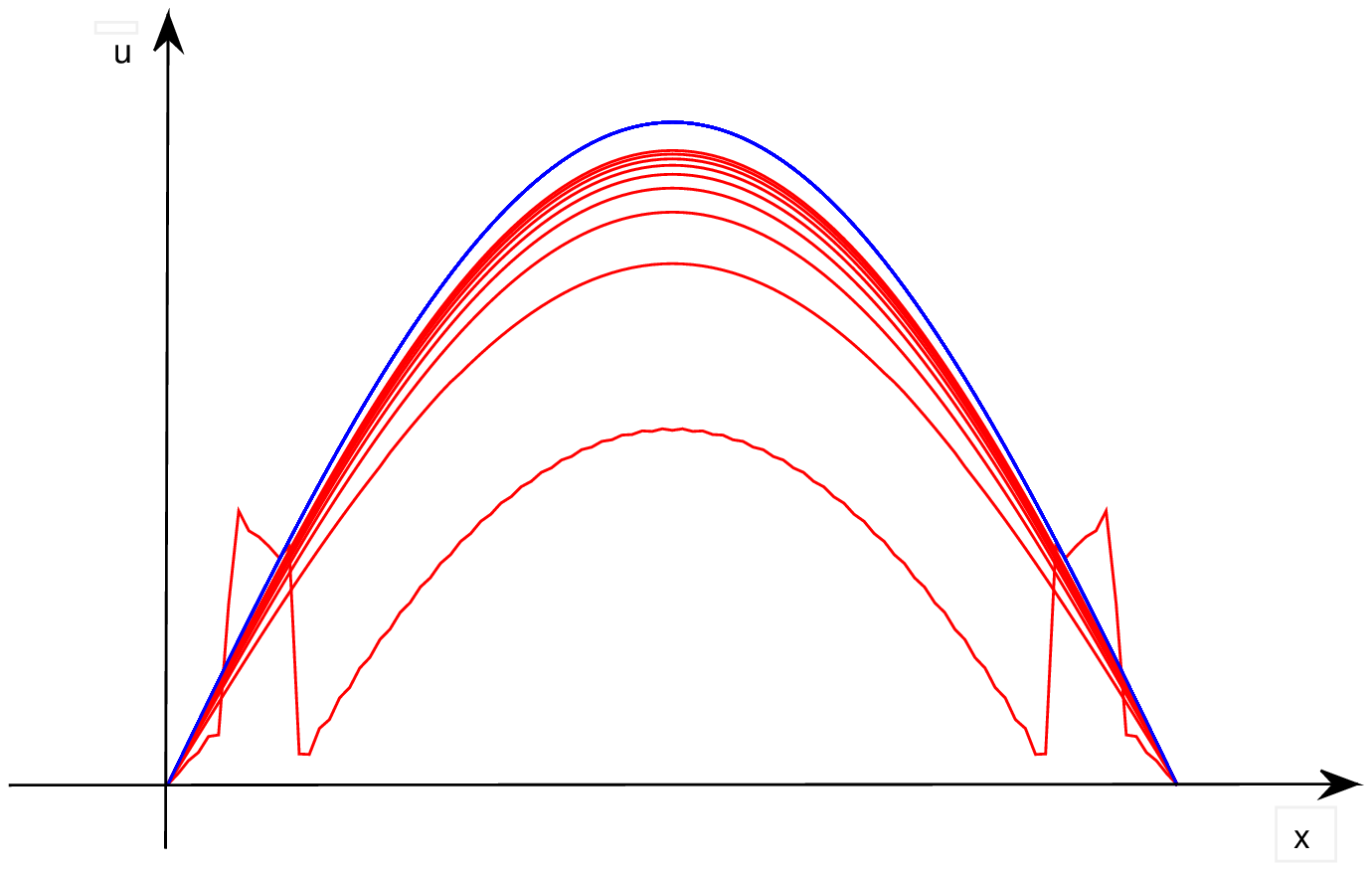}  
 \end{center}
\smallskip
\center{}
\end{figure}
 \end{center}

The previous plots suggest that the \GP equation   should approximate the wave equation. In fact, if $N(t)$ is ``constant enough'' we have a good approximation. To understand this point, we consider  \MC with 
$N(t)= e^{- {  k}t}$ (the initial condition is a step function centered at $t=\pi/2$, as in the previous figures). \Cref{Euris:4} on the left shows the plots of the solutions  at the time $T=2.5$, i.e. after the reflection, of the wave equation (blue)  and   \MC  (red) for various values of $k$, equal $4\,,  4/2\,,  \dots\,,  4/9$.

Even more interesting is the fact that  \GP can   be used also to approximate the heat equation (from the point of view of the plots of the solutions and not from the point of view of the analytic properties of the solutions) and this justifies the use of \GP as an hyperbolic model of heat diffusion proposed in~\cite{CattaneoMODENA1949,GurtinPipkinARCHIVE1968}. This is shown
by the plots in \cref{Euris:4} (right) which show the solution of the heat equation (at time $T=2.5$) and of \MC when the kernel is $N(t)={ k}e^{-{  k} t}$  $k=4\,,\ 8\,,\ 12\,,  \dots\,,\ 36$ (the initial condition is a step function centered at $t=\pi/2$, as in the previous figures).

  \section{The solution of the equation with memory}
 
  It is impossible and out of our goal to describe the many papers devoted to the study of the solutions of \cref{eq:memo} and similar equations. We confine ourselves to a short description.
  
  Due to the Volterra structure of the equation, a natural approach to the definition/analysis of the solutions is via Laplace transform. In fact, let the   initial   condition be
  \[
w(0)=w(x,0)=w_0(x)   \,.
  \]
  Then, the formal Laplace transform of \cref{eq:memo}  gives the following Laplace equation for $ \hat w(\zl)=\hat w(\zl,x)$, which has to be solved for every $\zl$ in a suitable right hand plane:
  \[
\zl \hat w( \zl)-\zaa A \hat w( \zl)-\hat N(\zl)A \hat w( \zl )=  w_0+\hat F( \zl)\qquad \hat w(\zl)=0\ {\rm on}\ \partial\ZOMq \,.
  \]
   If it is possible to prove that 
  \[ 
  \hat w(\zl)=\left (\zl I-\left (\zaa+\hat N(\zl)A\right )\right )^{-1}\left [w_0+\hat F(\zl)\right ]
  \]
  is the Laplace tranformation of a function $w( t)$, in a suitable space, then it is natural to choose such function $w$ as the solution of the problem. This idea has been used for example
   in~\cite{DaPratoPADOVA1980,DaPratoPADOVA1980-2,NarainJosephLINEARIZED1983} and in the important book~\cite{PruessLIBRO1993}. It is used in~\cite{DavisProcAMS1975,DavisSIAMAPPLMATH1976,BelleniBUMI1978}
to prove  finite propagation speed, when the kernel  has a strictly proper rational Laplace transform and in the already cited paper~\cite{RenardyREOLOGICAL1982}
 where, we recall, it is proved  that signals propagate with a finite speed  if $N(t)$ is smooth for $t\geq 0$ with $N(0)>0$ and also if if  $N(t)\sim 1/t^\zg$, $\zg\in(0,1)$ for $t\to 0^+$ but a (forward) wavefront exists only in the first case.
 
 See~\cite{CORDUlibro} for a different approach.

 When the boundary condition is nonhomogeneous, for example $ (w-Df)\in {\rm dom}\, A $ (for example when $A=\Delta$ with homogeneous Dirichlet boundary conditions  this amount to impose $ w(x,t)=f(x,t) $ $ x\in \partial\ZOMq $) then the solution is often defined by transposition, as for 
 example in 
  \cite{BenabABSTRACT2002,GamboaQAPPLMATH2016,KimSIAM1993,MunozADVDE2004}. 
 
When $ M(t) $ in~\cref{eq:memoSecondOrder} is of class $ H^2 $ (i.e. $ N\in H^3 $)
an important observation in~\cite{MacCamyMODEL1977,MacCamyINTEGRODIFFERENTIAL1977}, can be used  
to represent \cref{eq:memo} as\footnote{the function $ F(t) $ is not the same  and if the initial conditions $ w(0) =w_0$ and $ w'(0)=w_1 $ are not zero then $ F(t) $ depends also on the vectors $ w_0  $ and $ w_1 $.}
\begin{equation}\ZLA{eq:DopoMacCamy}
w''=A w+bw+\intt K(t-s) w(s)\ZD s +F(t)\,.
\end{equation} 
The trick is to formally solve the Volterra integral equation~(\ref{eq:memoSecondOrder}) respect to the ``unknown'' $ Aw $ followed by two integration by parts.  
This method, known as ``MacCamy trick'', is used to simplify the treatment of   \cref{eq:memo}  but (\ref{eq:memo})   can also  be studied directly, even without the help of MacCamy trick. We mention   an 
important paper,~\cite{BelleniBUMI1978}, which concerns \GP\!\!. In this paper,   the solutions of   \cref{eq:memo}  (with $ A=\Delta $ and homogeneous boundary conditions) are represented via   cosine operator theory.
   Using the cosine representation of the wave equation with boundary control in~\cite{LasieckaTriggianiCOSINE1981}, Belleni-Morante representation can be extended   to the case that the boundary condition is nonhomogeneous, as done in \cite{PandolfiAMO2005}.

   When $ \zaa>0 $ a similar idea, with $ e^{At} $ in the place of the cosine operator generated by $ A $, had already been used in~\cite{FriedmanTRANSAMS1967}.

Finally, we mention an important approach,  introduced in~\cite{DafermosJDE1970}, which reduces~(\ref{eq:memo}) to a semigroup system in a suitable ``state space''. This approach is of the utmost importance in the study of energies, the asymptotic properties, non linear perturbation and the definition of state (see~\cite{FabrizioARCH2010}) but it does not seem to be usable in the study of controllability.

In the special case that the kernel is a Prony sum, a a representation of the system which displays several ``components'' of the solution and which is related to the study of controllability is in~\cite{RenardySCL2005,ChowdhuryJMATHFLUIDMECH2017,WangGuoFU2009}.

\section{An overview of the results on controllability}

Interest on controllability of systems with memory arose in the '80-thies of the last century.
The oldest paper on controllability I know,~\cite{BaumeisterJMAA1983}, studies controllability for \CG\ but most of the significant results on controllability of \CG are recent while controllability of \GP has a longer and richer history. So, we consider \GP first. 
 
  The prototype of \GP (i.e. \cref{eq:memo} with $\zaa=0$ and $N(0)>0$) is often written in the form~(\ref{eq:memoSecondOrder})    or~(\ref{eq:DopoMacCamy}).
 
 We noted that, when $ K=0 $ and $A=\Delta$, \cref{eq:DopoMacCamy} reduces to the telegrapher's equation. We are not going to consider such special case but we mention the fact that  the telegrapher's equation  has the same controllability properties as the wave equation. In particular there exist  subsets of $ \partial \ZOMq $  and correspondigly times $ T $ such that the wave and the telegrapher's equation  are controllable in $ L^2(\ZOMq)\times  H^{-1} (\ZOMq) $ using square integrable controls applied to $ \Gamma $, see \cite{SorianoCOMPLUTENSE1995}, which gives a proof based on HUM method and, in the case ${\rm dim}\,\ZOMq=1$, the paper~\cite{ShubovSIAMcontr1997} which uses a moment operator approach.

  \subsection{The controllability of \GP equation}
The study of controllability   of \GP (not reduced to the telegrapher's equation) was initiated in  the PhD thesis of G. Leugering, supervised by W.~Krabs.  The first published paper 
is~\cite{LeugeringFADING1984}. This paper  studies controllability   in a $2$-dimensional region  $\ZOMq\subseteq\zzr^2$ (but the dimension $2$ seems then absorbed by the abstract formulation used in this paper)  under the action of an $L^2$ control \emph{everywhere and evenly  distributed in $\ZOMq$.} This means that the system is described by \cref{eq:memoSecondOrder} and $ F(t)=F(x,t)  \in L^2(0,T;L^2(\ZOMq) $ for every $ T>0 $ is the control. It is assumed that the operator $ -A $ is selfadjoint positive defined and it has an orthonormal sequence of eigenvectors $ \phi_n $ which span the space $ X=L^2(\ZOMq) $.
 
The paper is mostly devoted to derive a representation formula  for the solutions   under suitable assumptions (in particular, the resolvent kernel of $M(t)$  is square integrable and  completely monotonic).  

The special assumptions on the relaxation kernel and in particular the   the fact that the control is \emph{everywhere and evenly distributed in the region,}  have the following consequence:
let
$ u $
solve
\begin{equation}\ZLA{LaONDELeuge1}
u''=Au+F\qquad \mbox{$F$= control}
 \end{equation}
 (with zero initial conditions, like the solution $ w $).
The formula for the solutions easily shows that at any time  $ T>0 $ the difference $ \left (w^f(T),w'^f(T)\right )- 
\left (u^f(T),u'^f(T)\right )$ is quasinilpotent in $ X\times {\rm dom}\,(-A)^{1/2} $. It is clear that
\[ 
f\mapsto \left (u^f(T),u'^f(T)\right )
 \]
 is surjective in $ {\rm dom}\,(-A)^{1/2}\times X  $ and so
\[ 
f\mapsto \left (w^f(T),w'^f(T)\right )
 \]
 is surjective too.
 
 The proof of quasinilpotency uses in a crucial way the special fact that the control acts \emph{evenly and everywhere} on $ \ZOMq $. The proof in this paper cannot be adapted to more general cases, not even to the case 
    $ F=b(x)f(x,t) $ ($ f $ is the control and $ b(x) $ not a   constant); but, the idea seen in this special case has been used in almost every subsequent paper. The idea is as follows:
 we associate to (\ref{eq:memo}) (with $\zaa=0$)  or to its different versions~(\ref{eq:memoSecondOrder}) or~(\ref{eq:DopoMacCamy}) the \emph{associated ``wave'' equation}
 \begin{equation}
\ZLA{eq:associated}
u''=Au \qquad \mbox{with the same control as~\cref{eq:memo}} 
\end{equation}
 (in concrete cases  
$ A $ is a partial differential operator in a region $\ZOMq$ and the control may acts on a part of the  the boundary).

Let the initial conditions be zero and let $ u_f $, $ w_f $
be the solutions which correspond to the control $ f $. Let
\[ 
\Lambda_E f=\left (u_f(T),u_f'(T)\right )\,,\qquad 
\Lambda_V f=\left (w_f(T),w_f'(T)\right ) 
 \]
 (sometimes only the first components---the displacement\footnote{or the temperature or the concentration, according to the interpretation of the model}---is considered).
 
 The concrete cases that have been studied are cases in which it is already known that $ \Lambda_E $ is surjective and it is then proved that surjectivity is retained by $ \Lambda_V $ since 
 \begin{equation}
 \ZLA{eq:diffeDAFArePICCOLAcompatta}
  \Lambda_E-\Lambda_V 
  \end{equation}
   is ``small'' in a suitable sense. 
 We distinguish two cases
 \begin{itemize}
 \item[\Smal]  {\bf ``smallness'' in norm} since surjectivity is retained by perturbation of small norm;
 \item[\Comp]  ``smallness={\bf compactness}'' in the sense that $ \Lambda_V-\Lambda_E $ is {\bf compact.} Surjectivity is not retained under compact perturbation. It is only preserved the fact that the image is closed with finite codimension. An additional argument is then  needed to conclude surjectivity.
 \end{itemize}
 
A paper which do not use these ideas is~\cite{FuYongZhangJDE2009} (and also~\cite{FuZHANG-2,LeugeringLIQUIDM2AS1987,LeugeringBEAMAPPLANAL1986} discussed below). The higly technical paper~\cite{FuYongZhangJDE2009}   studies the case that $ A $ is a uniformly elliptic operator is a region $ \ZOMq $ and the control is a distributed control in $ \ZOMq $ (possibly localized close to a suitable part of the boundary). Controllability is proved via Carleman estimates of the adjoint system.

Both the approaches  have their merits:
\begin{itemize}
\item usually, the \Smal  approach require that the memory kernel have a ``small'' norm, but it may be possible to prove controllability for kernels which have low regularity;
\item the \Comp approach can be used for large memories, but more regularity is required.
\end{itemize}
 
 The papers which uses the approach \Smal assume that the memory kernel is ``small''. Often it is assumed that $ N(t)=  \ZEP N_0(t)$  and typical statements are the existence of $ \ZEP_0 $ such that controllability holds for $ 0\leq \ZEP<\ZEP_0 $ as in the final chapter of~\cite{LagneseLionsLIBRO1988}, which studies controllability of a Kirchhoff plate with memory when the control acts in a particular boundary conditions.  Also the paper~\cite{LasieckaKirchhoffCONF1989} studies controllability of a Kirchhoff plate with memory when the control acts in a more natural  boundary conditions. The special feature of the model in~\cite{LasieckaKirchhoffCONF1989} is that the operator~(\ref{eq:diffeDAFArePICCOLAcompatta}) is small in norm (if $\ZEP<\ZEP_0$) but not compact.
 
 Of course $\ZEP_0$ depends on $N_0(t)$. The dependence is made explicitly in 
 the paper\footnote{see also~\cite{LoretiCOMPRENDMATH2009,Loreti2011261}.}~\cite{LoretiJDE2010} which
 studies controllability of \MC\!\!. It is proved that controllability holds provided that $\beta<\ZEP_0=\frac{2}{3}\eta$ thanks to delicate estimates and the use of Ingham and Haraux type inequalities.

The paper~\cite{CavalcantiELECTRJ1998} proves controllability of a viscoelstic incompressible body  
 under a ``smallness'' assumption of the memory kernel and of its derivative.

  \Smal has been used  also in the study of interconnected systems, when one of the connected systems has a memory, as in~\cite{MunozADVDE2003,MunozADVDE2004}. The paper~\cite{GamboaQAPPLMATH2016} studies an interconnected system and the ``smallness'' assumption is not on the memory kernels but on the ``strength'' of the interconnection.
 
 In the special case that $ Au=u_{xx} $ or $Au=u_{xxxx}$ (viscoelastic string or beams)
it is also possible to study controllability directly, without relaying on either \Smal or \Comp\!\!, as in~\cite{LeugeringLIQUIDM2AS1987}. This paper  reduces    controllability (of a viscoelastic string)   to the solution of a Volterra integral equation of the second kind thanks to D'Alembert formula. A similar approach is in~\cite{BarbuDIFFINTEQ2000} whose final section\footnote{the main part of this paper is concerned with \CG as discussed below.} studies reachability of smooth targets when the  controls are evenly distributed on a subset of the interval occupied by a viscoelastic string.
In~\cite{LeugeringBEAMAPPLANAL1986}  controllability is proved via abelian theorems for the  Laplace transform. This paper  studies  controllability under boundary control, as~\cite{LeugeringDYNAstabSYSTEMS1989,LeugeringINTEGRODIFFERENTIAL1987}. More specifically, the paper~\cite{LeugeringBEAMAPPLANAL1986} studies controllability of a viscoelastic \emph{beam} under the action of a boundary control, i.e. the equation is now
\[
w_{tt}=w_{xxxx} +\intt M(t-s) w_{xxx}(s)\ZD s
\qquad\left\{\begin{array}{l}
w(t,0)=w(t,l)=w_{xx}(t,0)=0\,,\\   w_{xx}(t,l)=f(t)=\mbox{control}\,.
\end{array}\right.
\]
Signals propagate with infinite speed, a property which permits controllability at an arbitrary small time (as first proved in~\cite{ZuazuaTEMPSPETITparis1987} for systems without memory). Infinite velocity of propagation is 
  inherited by the viscoelastic system and exact controllability is proved in~\cite{LeugeringBEAMAPPLANAL1986}, provided that $ M\in W^{1,2}(0,T)$   for every $T$ and $M(0)=0$, in a suitable space defined in terms of the spectral properties of the bilaplacian. The proof of controllability rests on estimates of the eigenfunction of the bilaplacian derived in~\cite{KrabsLeugeringSeidman1985}. These estimates are used in order to prove that the adjoint of the controllability operator is coercive\footnote{``inverse inequality'' is a term often used to denote coercivity of the adjoint of the control map.}. The proof of this fact  is via the Laplace transform of the adjoint equation. The Laplace transform of the solutions   of the adjoint equation solve   
  \begin{equation}
  \ZLA{eq:elliEqueBEAM}
 \zl^2\hat w_{\zl, xxx}=-\frac{\zl^2}{1+\hat M(\zl)}\hat w(\zl, x)+\frac{1}{1+\hat M(\zl)} H(s)
  \end{equation}
   at every complex frequency $\lambda$ ($H$ computed from the associated equation with the memory $M$ put equal zero). The assumption that ${\rm dim }\,\ZOMq=1$ is used in order to explicitly solve this equation for every $\zl$ (in a right half plane).

The analogous of \cref{eq:elliEqueBEAM} can be explicitly solved also in domains $\ZOMq$  with ${\rm dim}\,\ZOMq>1$, provided that $\ZOMq$ has a  special form. In particular if $\ZOMq$ is a square  $\ZOMq=(0,l)\times (0,l)$. This case is treated with similar methods in~\cite{LeugeringINTEGRODIFFERENTIAL1987}. The case of general (smooth) regions with one boundary control in different boundary conditions is   in~\cite{KIMplate1992,PandolfiSIAM2017} via \Comp (the paper~\cite{KIMplate1992} is discussed below).

Finally we cite~\cite{LeugeringONEDIMTHERMOVISCO1985,LeugeringM2AS1986,LeugeringLIQUIDM2AS1987}.
  Noticeable following papers   study controllability   in a (smooth) region $\ZOMq\subseteq\zzr^n$, with boundary control in the deformation (i.e. in the Dirichlet boundary condition). In particular,~\cite{YanCOMPLUTENSE1992} and~\cite{KimSIAM1993} appeared in the early 90-thies of the last century (the method in this paper are extended to the control of a plate equation in~\cite{KIMplate1992}). The \Comp approach is used in these papers. The paper~\cite{YanCOMPLUTENSE1992} studies 
  system (\ref{eq:DopoMacCamy}) with $b= 0$ and smooth kernel\footnote{in fact, in the papers~\cite{YanCOMPLUTENSE1992,KIMplate1992,KimSIAM1993}  the kernel is $K(t,s)$ instead of $K(t-s)$ and $b=b(t)$ in the papers~\cite{KIMplate1992,KimSIAM1993}    but this  facts do not add difficulties to the problem. The paper~\cite{YanCOMPLUTENSE1992} implicitly assume the MacCamy transformation already done.} (of class $W^{1,1}(0,T)$).  Both these papers lift the
  \emph{direct inequality}\footnote{i.e. the fact that the trace on $\partial\ZOMq$ of the normal derivative of a wave equation (with homogeneous boundary conditions and initial data in $H^1_0(\ZOMq)\times L^2(\ZOMq)$) is square integrable.}   from the wave equation to the \GP equation.

As we said,~\cite{YanCOMPLUTENSE1992} assume more regularity of the kerne, while~\cite{KIMplate1992,KimSIAM1993} assume that the control act on the entire boundary.
 In fact the the direct inequality holds without this regularity assumption\footnote{the direct inequality has been proved for \MC in~\cite{LoretiPARMA2016} using multiplier methods. As ``admissibility'' of the input or of the output operator, it has been studied for example in~\cite{JacobJMAA2007,HaakJDE2009,JacobBirgitJFA2010,JungSIAMJCONTR2000} by using techniques of harmonic analysis. An Hautus' type condition, as examined in~\cite{DAPRATOJINTEQAPPL1990}, should have a role in this approach.}, as proved in~\cite{KIMplate1992,KimSIAM1993,PandolfiJINTEQAPPL2015,PandolfiPARMA2016,PandolfiSIAM2017,PandolfiESAIM2017}.
  After that, these references use 
  a compactness/uniqueness argument    to derive controllability.  
  
  As we noted,  the paper~\cite{KIMplate1992,KimSIAM1993} base the proof on the fact that    the control acts on the entire boundary while~\cite{PandolfiJINTEQAPPL2015,PandolfiPARMA2016,PandolfiSIAM2017,PandolfiESAIM2017} use  solely controllability of the associated ``wave'' equation  without assumptions on the active part of the boundary (a part the geometric control condition  in~\cite{BardosSIAMJCONT1992,BurqASYMPANAL1997}). The paper~\cite{PandolfiSIAM2017,PandolfiESAIM2017} proves controllability of the plate and Navier equations.
   
   The reason for this assumption  that the control acts on the entire boundary in~\cite{KIMplate1992,KimSIAM1993} is the use of a    unique continuation result for the adjoint equation, which uses a delicate analysis of the extension of a solution which is zero on the boundary, with normal derivative equal zero too, to a larger region.

My interest on this subject was stimulated by Maria Grazia Naso talk at the IFIP-TC7 Conference in Antibes (2003), where she presented the results in~\cite{MunozADVDE2003}. The paper~\cite{PandolfiAMO2005}  was inspired by this talk. This paper in particular extends Belleni-Morante approach (in~\cite{BelleniBUMI1978}) to the representation of the solutions of \GP with  boundary controls  and it uses the obtained representation formula in order to prove compactness of the operator~(\ref{eq:diffeDAFArePICCOLAcompatta}). Surjectivity is then proved via   Laplace transformation of the solutions and Baire Theorem. The ``sharp control time'' was not   identified in this paper (see a correction to the last part of the paper in~\cite{PandolfiErratumAMO2011}).

A different use of compactness is in the papers~\cite{PandolfiINTEQOPTH2009,PandolfiDCDSSerB2010}\footnote{see also~\cite{LoretiSIAM2012,RomanovShamev}.}. these papers consider controllability of a string  with boundary control and identify a special sequence related to controllability\footnote{these sequences are obtained by first projecting the solutions on the Dirichlet eigenfunctions of the laplacian,  as in~\cite{LeugeringSTRING1989}.}. Controllability holds if and only if this sequence is a Riesz sequence, i.e. if and only if it can be transformed to an orthonormal sequence by a linear bounded boundedly invertible transformation.  This approach has then been extended in~\cite{PandolfiJINTEQAPPL2015}  to bodies occupying bounded regions regions in $\zzr^n$ (with smooth boundary. In the applications of course $n\leq 3$).

One of the sequences encountered in this papers (for bodies in regions of $\zzr^n$) is $ \{z_n(t)\zg_1\phi_n \} $ ($\zg_1$ is the trace of the normal derivative) where $z_n(t)$ solves 
 
\[ 
z_n'(t)=-\zl_n^2\int_0^t N(t-s) z_n(s)\ZD s\,,\qquad z_n(0)=1\,.
 \]

The fact that  $ \{z_n(t)\zg_1\phi_n \} $ is a Riesz sequence is a consequence of the controllability of the associated wave equation, as proved via  Bari Theorem, hence   using a compactness argument (in fact, it is true that the full force of the Bari Theorem is not needed in the proof, see~\cite[Th.~3.10 and Ex.~3.15]{PandolfiLIBRO2014}).

In~\cite{PandolfiJINTEQAPPL2015}, this Riesz sequence has been   used to identify the infimum of the control times.

The operator approach in~\cite{PandolfiAMO2005} and the moment approach in~\cite{PandolfiJINTEQAPPL2015}   have been merged to prove controllability in three dimensional viscoelasticity and for viscoelastic plates in \cite{PandolfiSIAM2017,PandolfiESAIM2017}.

Finally, we note that the papers we cited up to now essentially study controllability when the control is either distributed within the region occupied by the body, or it affects the boundary deformation.  Boundary control in the Dirichlet condition is important for the solution of inverse problems (we are not going to review these results). Applications to mechanics are more likely to control via the boundary traction. Control under   boundary traction of a viscoelastic  string   is studied in~\cite{NegrescuAPPLIED2016}.

\subsection*{Controllability of the stress or flux}

The problems considered in the papers cited up to now concern controllability of   the pair $(w(T),w'(T))$ (or of solely $w(T)$). 
According to the different applications, we studied controllability of (deformation/velocity of   deformation), controllability of the temperature or of the concentration.

As noted in~\cite{DubovaCRPARIS1998,RenardySCL2005,Renardy2008}, controllability of the traction  is   more important then the control of velocity, since for example traction can survive   solidification of a body and it can be the cause of fractures.

When $A=\Delta$ the traction is 
\[ 
\ZSI(x,t)=-\left [N(0)\nabla w(x,t)+\intt N'(t-s)\nabla w(x,s)\ZD s\right ]
 \]
 while in thermodynamics
 \[ 
 q(x,t)=-\intt N(t-s)\nabla w(x,s)\ZD s
  \]
  is the flux.

When $x\in (a,b)$ and for general kernels, controllability of the pair temperature/flux has been proved in~\cite{AvdoninQAPLLMATH2013}.
and interpreted as independence of temperature and flux in~\cite{Avdonin2012} Note that in the standard heat equation, i.e. when Fourier law is assumed,  the flux is determined by the temperature. This controllability result can be interpreted as follows (see~\cite{Avdonin2012}):  temperature and flux are independent for \GP  when $T$ is   large.

Controllability of the pair (stress,deformation)  and (stress,velocity) is studied in~\cite{PandolfiEVOLEQCONTTH2013} where it is proved that the pair (stress,velocity)   is controllable after a suitable time determined by the velocity of wave propagation.   Hence,  stress and velocity of deformation are independent  after that a sufficiently large time is elapsed.  Instead,   stress and deformation are not independent, but the stress and an arbitrary number of the component of the deformation (respect to the basis of the eigenfunctions) can be arbitrarily assigned. Again we note that when Hooke law hold, deformation uniquely identifies the resulting stress at every instant of time.

The arguments in~\cite{AvdoninQAPLLMATH2013,PandolfiEVOLEQCONTTH2013}
are based on moment methods and Riesz sequences, via Bari Theorem; hence they fall in the approach \Comp to controllability. 

In the contest of the \J equation, this problem has been studied for example in~\cite{BoldriniSIAM2012, DubovaSCL2012}. The papers \cite{ChowdhuryJMATHFLUIDMECH2017,Mitra2018575,RenardySCL2005} study the case that the kernel is a Prony sum, and study controllability of several ``components'' of the stress.

  \subsection{Controllability of the \CG equation}

  As stated already, the first paper on controllability of the \CG equation is \cite{BaumeisterJMAA1983} which proves \emph{lack of exact controllability.}
    The next significant paper,~\cite{BarbuDIFFINTEQ2000}, appeared almost thirty years later.  
  This paper considers the case $ A=\Delta  $ in a region $ \ZOMq\subseteq \zzr^n $ with distributed control acting in a subregion 
  $ \ZOM\subseteq \ZOMq $, i.e. 
   \[ 
 F(x,t)=  {\bf  1}_{_\ZOM}(x) f(x,t)\quad \ZOM\subseteq \ZOMq\qquad \mbox{$f(x,t)$= control}\,.
  \]
  This paper proves approximate controllability in $ L^2(\ZOMq) $ by using Laplace transform and properties of analytic functions when
  \[ 
  N(t)=\sum _{j=1}^n  a_j e^{-\zaa_j t} +\sum _{j=1}^m \int_{I_j} b_j(s) e^{-st }\ZD s
   \]
 where $ \zaa_j>0 $, $a_j\geq 0$;  $ I_j $ are bounded intervals and $ b_k(s)>c_k>0 $ on $ I_k $ is integrable\footnote{approximate controllability when the relaxation kernel has compact support is in~\cite{LavanyaBALACH}.}. 
 
  In fact, approximate controllability in $ L^2(\ZOMq) $ holds for every   $ H^1 $ relaxation kernel, as proved in~\cite{HalanayJMAA2015}.

The solutions of \cref{eq:memo} with $ \zaa>0 $
 are given by
 \begin{equation}\ZLA{Eq:voltePERCGdaINTEGR} 
w(t)=u(t)+\intt e^{\zaa A (t-s)}\int_0^s N(s-r) A w(r)\ZD r\,\ZD s
  \end{equation}
  where 
  $ u(t) $ solves the ``heat'' equation~(\ref{eq:memo}) with $ N=0 $, with the same initial condition and control as \cref{eq:memo}. An integration by parts in \cref{Eq:voltePERCGdaINTEGR}  gives the following Volterra integral equation for $ w(t) $
 
  \begin{align}
\nonumber  w(t)=u(t)-\frac{1}{\zaa} \intt N(t-r) w(r)\ZD r\\
\nonumber   +\frac{1}{\zaa}\intt e^{\zaa A(t-s)}\left [
N(0) w(s)+\ints N'(s-r) w(r)\ZD r  
  \right ]\ZD s=\\
\nonumber=  \underbrace{\left [e^{At} \xi-A\intt e^{A(t-s) }Df(s)\ZD s\right ]}_{u(t)}
 -\frac{1}{\zaa} \intt N(t-r) w(r)\ZD r\\
 +\frac{1}{\zaa}\intt e^{  A(t-s)}\left [
N(0) w(s)+\ints N'(s-r) w(r)\ZD r    \right ]\ZD s
  \ZLA{Eq:voltePERCGintegrata}
  \end{align}

   It is clear  that $t\mapsto w(t)$ has the same ``regularity''  properties  as the solution of the heat equation. \emph{In particular, if $ A=\Delta $ on a region $ \ZOMq $ and the square integrable control acts on the boundary, $  w(t) $ is only square integrable, but it is not continuous.} So, ``reachability'' at time $ T $ makes sense only if the control $ f(t) $ has the property that $ w(t) $ is continuous 
  at $ T $ (from the left). It is proved in~\cite{HalanayJMAA2015,HalanayDCDS-A2015} that this class of control is sufficiently rich so to have approximate controllability  in $ L^2(\ZOMq) $ for every $ T>0 $ (see also~\cite{HalanayDCDS-A2015}). It is also possible to see from the representation~(\ref{Eq:voltePERCGintegrata})   that ``irregular'' targets cannot be reached, i.e. exact controllability is impossible (as already proved in~\cite{BaumeisterJMAA1983}) precisely  as in the case of the heat equation.
  
  The heat equation     is also controllable to zero: fix  any final time  time $ T>0 $. There exists a control $ f $ (either distributed or acting on a part of the boundary) such that the solution  $ u $ of the heath equation, i.e. the equation with $ N=0 $ is continuous at $ T $ and satisfies $ u(T)=0 $ (and $T>0$ can be taken arbitrarily small).
  
 It is  natural to conjecture   that this property is inherited by the system with memory, i.e. that it is possible to drive $w(t)$ to hit the target $0$, $w(T)=0$, but in fact it is not so. 
This conjecture was disproved by exhibiting   suitable counterexamples in~\cite{GuerreroESAIM2013,HalanaySCL2012} and then for every $H^1$ kernel in~\cite{HalanayJMAA2014,HalanayJMAA2015}.

   \section{Controllability with moving controls}
   
   We cited the state space approach introduced by Dafermos in~\cite{DafermosJDE1970} in order to have a ``time invariant'' system, which can be described via semigroup theory. In this approach the equation   take into account the entire past history, from 
   the ``creation of the universe'' , i.e. \cref{eq:memo} is  rewritten in the form
    \begin{equation} \ZLA{eq:memoTOTALE}
w'=\alpha A w(s)+\int_{-\ZIN}^t N(t-s) A w(s)\ZD s  +F(t) 
 \end{equation}
 and the state at time $ t $ has a component which is the entire history $ w(s) $, $ s<T $. It is clear that state space exact controllability is impossible due to causality, since the control cannot affect the past evolution, before it was applied. So, the kind of controllability we examined so far is not exact controllability in the state space, but it is a kind of \emph{relative controllability}, in a sense once used in the theory of systems with delays in $ \zzr^n $.  In spite of this limitation, it turns out that this kind of \emph{relative controllability } is usefull to solve certain inverse problems (see~\cite{PandolfiDCDSserS2011,PandolfiANNREVCONT2010,PandolfiJMAA2013}, see~\cite{KabanikhinLorenziBOOK1999} for different approaches). But, it is clear that in certain applications, for example in the theory of the quadratic regulator problem on a half line, the following property would be important: \emph{The initial vector $ \xi  $
 is {\em controllable to the rest\/} if we can find $ T>0 $ and a control $ u(t) $ with $ u(t)=0 $ for $t>T$
such that the corresponding solution $ w(t) $ satisfies $ w(t)=0 $ for $ t>T $.
{\em The system is   controllable to the rest\/} when every initial vector $ \xi $ is controllable to the rest.
 }

 Clearly, controllability to the rest is achievable only in extremely special cases (in fact, essentially when the equation is the integrated version of a wave or telegrapher's equation, see~\cite{IvanovJMAA2009}). From the other side,   certain system without memory  which are approximately controllable may  lack controllability to zero, hence also controllability to the rest.
 An example is in~\cite{MartinRosierSIAM2013} where a cure to this fact is provided:  the system is controllable to zero provided that the control is applied in a region inside the body, which is not fixed but which moves in a suitable way (see also~\cite{ChavezRosierLMPA2014}).
 Stimulated by these results,   systems with persistent memory have been examined from the point of view of controllability with moving controls in~\cite{ChavezSIAM2017,LuJMATHPURAPPL2017}.
In these papers the two models \CG and \GP are studied, with homogeneous Dirichlet conditions and
\[
F(x,t)={\bf 1}_{\omega(t)} g(x,t)
\]
where $g$ is the control, acting on the moving region $\ZOM(t)\subseteq \ZOMq$.
   It is proved that if the set $ \ZOM(t) $ moves inside the set $ \ZOMq $ and spans all of $ \ZOMq $ in a suitable way then it is possible to achieve controllability to zero of the quantities:
   \begin{itemize}  
  \item[~]  \colorbox{white}{\GP model:}  displacement, velocity   and (integrated) stress:  
  \end{itemize}
   
  \[
\left (w(T), w'(T),\int_0^T N(T-s)\nabla w(s)\ZD s\right ) 
  \]
 \begin{itemize}
  \item[~] \colorbox{white}{\CG model:}   temperature and flux i.e.
\end{itemize}

 \[
\left (w(T),\int_0^T N(T-s)\nabla w(s)\ZD s\right )\,. 
  \]

   This null controllability properties imply controllability to rest if  
   \[ 
   w(T)=0\,,\quad \int_0^T N(T-s)w(s)\ZD s=0\ \implies\ \mbox{$\intt N(t-s) w(s)\ZD s=0$ for $t>T$} 
    \]
  as it is the case when $N(t)=e^{-\beta t}$.

  \bibliographystyle{plain}
 \bibliography{bibliomemoria}{ }

\begin{thebibliography}{100}

\bibitem{Amendola2012}
G.~Amendola, M.~Fabrizio, and J.M. Golden.
\newblock {\em Thermodynamics of materials with memory: Theory and
  applications}.
\newblock {S}pringer, New York, 2012.

\bibitem{Avdonin2012}
S.~Avdonin and L.~Pandolfi.
\newblock Temperature and heat flux dependence/independence for heat equations
  with memory.
\newblock {\em Lecture Notes in Control and Information Sciences}, 423:87--101,
  2012.

\bibitem{AvdoninQAPLLMATH2013}
S.~Avdonin and L.~Pandolfi.
\newblock Simultaneous temperature and flux controllability for heat equations
  with memory.
\newblock {\em Quarterly of Applied Mathematics}, 71(2):339--368, 2013.

\bibitem{BarbuDIFFINTEQ2000}
V.~Barbu and M.~Iannelli.
\newblock Controllability of the heat equation with memory.
\newblock {\em Differential Integral Equations}, 13(10-12):1393--1412, 2000.

\bibitem{BardosSIAMJCONT1992}
C.~Bardos, G.~Lebeau, and J.~Rauch.
\newblock Sharp sufficient conditions for the observation, control, and
  stabilization of waves from the boundary.
\newblock {\em SIAM J. Control Optim.}, 30(5):1024--1065, 1992.

\bibitem{BaumeisterJMAA1983}
J.~Baumeister.
\newblock Boundary control of an integro-differential equation.
\newblock {\em J. Math. Anal. Appl.}, 93(2):550--570, 1983.

\bibitem{BelleniBUMI1978}
A.~Belleni-Morante.
\newblock An integro-differential equation arising from the theory of heat
  conduction in rigid materials with memory.
\newblock {\em Boll. Un. Mat. Ital. B (5)}, 15(2):470--482, 1978.

\bibitem{BenabABSTRACT2002}
A.~Benabdallah and M.~G. Naso.
\newblock Null controllability of a thermoelastic plate.
\newblock {\em Abstr. Appl. Anal.}, 7(11):585--599, 2002.

\bibitem{BoldriniSIAM2012}
J.~L. Boldrini, A.~Doubova, E.~Fern\'andez-Cara, and M.~Gonz\'alez-Burgos.
\newblock Some controllability results for linear viscoelastic fluids.
\newblock {\em SIAM J. Control Optim.}, 50(2):900--924, 2012.

\bibitem{BoltzmannWIN1874}
L.~Boltzmann.
\newblock Zur theorie der elastischen nachwirkung.
\newblock {\em Wien. Ber.}, 70:275--306, 1874.

\bibitem{BoltzmannWIED1878}
L.~Boltzmann.
\newblock Zur theorie der elastischen nachwirkung.
\newblock {\em Wied. Ann.}, 5:430--432, 1878.

\bibitem{BonaccorsiJEVOLEQ2012}
S.~Bonaccorsi and G.~Desch.
\newblock Controllability of a class of {V}olterra equations in {H}ilbert
  spaces with completely monotone kernel.
\newblock {\em J. Evol. Equ.}, 12(3):675--695, 2012.

\bibitem{BucciPandolfiArxive}
F.~{Bucci} and L.~{Pandolfi}.
\newblock {On the regularity of solutions to the Moore-Gibson-Thompson
  equation: a perspective via wave equations with memory}.
\newblock {\em ArXiv e-prints}, December 2017.

\bibitem{BurqASYMPANAL1997}
N.~Burq.
\newblock Contr\^olabilit\'e exacte des ondes dans des ouverts peu r\'eguliers.
\newblock {\em Asymptot. Anal.}, 14(2):157--191, 1997.

\bibitem{CattaneoMODENA1949}
C.~Cattaneo.
\newblock Sulla conduzione del calore.
\newblock {\em Atti Sem. Mat. Fis. Univ. Modena}, 3:83--101, 1949.

\bibitem{CavalcantiELECTRJ1998}
M.~M. Cavalcanti, V.~N. Domingos~Cavalcanti, A.~Rocha, and J.~A. Soriano.
\newblock Exact controllability of a second-order integro-differential equation
  with a pressure term.
\newblock {\em Electron. J. Qual. Theory Differ. Equ.}, pages No. 9, 18, 1998.

\bibitem{ChavezRosierLMPA2014}
F.~W. Chaves-Silva, L.~Rosier, and E.~Zuazua.
\newblock Null controllability of a system of viscoelasticity with a moving
  control.
\newblock {\em J. Math. Pures Appl. (9)}, 101(2):198--222, 2014.

\bibitem{ChavezSIAM2017}
F.~W. Chaves-Silva, Xu~Zhang, and E.~Zuazua.
\newblock Controllability of evolution equations with memory.
\newblock {\em SIAM J. Control Optim.}, 55(4):2437--2459, 2017.

\bibitem{ChowdhuryJMATHFLUIDMECH2017}
S.~Chowdhury, D.~Mitra, M.~Ramaswamy, and M.~Renardy.
\newblock Approximate controllability results for linear viscoelastic flows.
\newblock {\em Journal of Mathematical Fluid Mechanics}, 19(3):529--549, 2017.

\bibitem{ChristensenBOOK1982}
R.M. Christensen.
\newblock {\em Theory of viscoelasticity, an introduction}.
\newblock Academic Press, New York, 1982.

\bibitem{ColemanGurtinZEITSCH1967}
B.~D. Coleman and M.~E. Gurtin.
\newblock Equipresence and constitutive equations for rigid heat conductors.
\newblock {\em Z. Angew. Math. Phys.}, 18:199--208, 1967.

\bibitem{CORDUlibro}
C.~Corduneanu.
\newblock {\em Integral equations and applications}.
\newblock Cambridge University Press, Cambridge, 1991.

\bibitem{CurieANNALES1889}
M.J. Curie.
\newblock Recherches sur la conductibilit\'e des corps cristallines.
\newblock {\em Ann. Chim. Phys.}, 18:203--269, 1889.

\bibitem{DaPratoPADOVA1980-2}
G.~Da~Prato and M.~Iannelli.
\newblock Linear abstract integro-differential equations of hyperbolic type in
  {H}ilbert spaces.
\newblock {\em Rend. Sem. Mat. Univ. Padova}, 62:191--206, 1980.

\bibitem{DaPratoPADOVA1980}
G.~Da~Prato and M.~Iannelli.
\newblock Linear integro-differential equations in {B}anach spaces.
\newblock {\em Rend. Sem. Mat. Univ. Padova}, 62:207--219, 1980.

\bibitem{DAPRATOJINTEQAPPL1990}
G.~Da~Prato and A.~Lunardi.
\newblock Stabilizability of integrodifferential parabolic equations.
\newblock {\em J. Integral Equations Appl.}, 2(2):281--304, 1990.

\bibitem{DafermosJDE1970}
C.~M. Dafermos.
\newblock An abstract {V}olterra equation with applications to linear
  viscoelasticity.
\newblock {\em J. Differential Equations}, 7:554--569, 1970.

\bibitem{DavisProcAMS1975}
P.~L. Davis.
\newblock Hyperbolic integrodifferential equations.
\newblock {\em Proc. Amer. Math. Soc.}, 47:155--160, 1975.

\bibitem{DavisSIAMAPPLMATH1976}
P.~L. Davis.
\newblock On the hyperbolicity of the equations of the linear theory of heat
  conduction for materials with memory.
\newblock {\em SIAM J. Appl. Math.}, 30(1):75--80, 1976.

\bibitem{DEkEELIUHinestroza}
D.~De~Kee, Q.~Liu, and J.~Hinestroza.
\newblock ({N}on-fickian) diffusion.
\newblock {\em The Canada J. of chemical engineering}, 83:913--929, 2005.

\bibitem{DelloroAMO2017}
F.~Dell'Oro and V.~Pata.
\newblock On the {M}oore-{G}ibson-{T}hompson equation and its relation to
  linear viscoelasticity.
\newblock {\em Appl. Math. Optim.}, 76(3):641--655, 2017.

\bibitem{DeschINTEGRALEQ1985}
W.~Desch and R.~Grimmer.
\newblock Invariance and wave propagation for nonlinear integro-differential
  equations in {B}anach spaces.
\newblock {\em J. Integral Equations}, 8(2):137--164, 1985.

\bibitem{DubovaSCL2012}
A.~Doubova and E.~Fern\'andez-Cara.
\newblock On the control of viscoelastic {J}effreys fluids.
\newblock {\em Systems Control Lett.}, 61(4):573--579, 2012.

\bibitem{DubovaCRPARIS1998}
A.~Dubova, E.~Fern\'andez-Cara, and M.~Gonz\'alez-Burgos.
\newblock Controllability results for discontinuous semilinear parabolic
  partial differential equations.
\newblock {\em C. R. Acad. Sci. Paris S\'er. I Math.}, 326(12):1391--1395,
  1998.

\bibitem{FabrizioARCH2010}
M.~Fabrizio, C.~Giorgi, and V.~Pata.
\newblock A new approach to equations with memory.
\newblock {\em Arch. Ration. Mech. Anal.}, 198(1):189--232, 2010.

\bibitem{FabrizioMorroCBMS}
M.~Fabrizio and A.~Morro.
\newblock {\em Mathematical problems in linear viscoelasticity}, volume~12 of
  {\em SIAM Studies in Applied Mathematics}.
\newblock Society for Industrial and Applied Mathematics (SIAM), Philadelphia,
  PA, 1992.

\bibitem{FattoriniBOOKcauchy1983}
H.~O. Fattorini.
\newblock {\em The {C}auchy problem}, volume~18 of {\em Encyclopedia of
  Mathematics and its Applications}.
\newblock Addison-Wesley Publishing Co., Reading, Mass., 1983.

\bibitem{FisherGurtin1965}
G.~M.~C. Fisher and M.~E. Gurtin.
\newblock Wave propagation in the linear theory of viscoelasticity.
\newblock {\em Quart. Appl. Math.}, 23:257--263, 1965.

\bibitem{FriedmanTRANSAMS1967}
A.~Friedman and M.~Shinbrot.
\newblock Volterra integral equations in {B}anach space.
\newblock {\em Trans. Amer. Math. Soc.}, 126:131--179, 1967.

\bibitem{FuYongZhangJDE2009}
X.~Fu, J.~Yong, and X.~Zhang.
\newblock Controllability and observability of a heat equation with hyperbolic
  memory kernel.
\newblock {\em J. Differential Equations}, 247(8):2395--2439, 2009.

\bibitem{GamboaQAPPLMATH2016}
P.~Gamboa, V.~Komornik, and O.~Vera.
\newblock Partial reachability of a thermoelastic plate with memory.
\newblock {\em Quarterly of Applied Mathematics}, 74(2):235--243, 2016.

\bibitem{GraffiCIMENTO1928}
D.~Graffi.
\newblock Sui problemi della eredit\`a lineare.
\newblock {\em Nuovo Cimento}, 5:53--71, 1928.

\bibitem{GraffiLOMBARDO1936}
D.~Graffi.
\newblock Sopra alcuni fenomeni ereditari dell'elettrologia.
\newblock {\em Rend. Istit. Lombardo Sci. Lett}, 68-69:124--139, 1936.

\bibitem{GuerreroESAIM2013}
S.~Guerrero and O.~Y. Imanuvilov.
\newblock Remarks on non controllability of the heat equation with memory.
\newblock {\em ESAIM Control Optim. Calc. Var.}, 19(1):288--300, 2013.

\bibitem{GurtinPipkinARCHIVE1968}
M.~E. Gurtin and A.~C. Pipkin.
\newblock A general theory of heat conduction with finite wave speeds.
\newblock {\em Arch. Rational Mech. Anal.}, 31(2):113--126, 1968.

\bibitem{HaakJDE2009}
B.~H. Haak, B.~Jacob, J.~R. Partington, and S.~Pott.
\newblock Admissibility and controllability of diagonal {V}olterra equations
  with scalar inputs.
\newblock {\em J. Differential Equations}, 246(11):4423--4440, 2009.

\bibitem{HalanaySCL2012}
A.~Halanay and L.~Pandolfi.
\newblock Lack of controllability of the heat equation with memory.
\newblock {\em Systems Control Lett.}, 61(10):999--1002, 2012.

\bibitem{HalanayJMAA2014}
A.~Halanay and L.~Pandolfi.
\newblock Lack of controllability of thermal systems with memory.
\newblock {\em Evol. Equ. Control Theory}, 3(3):485--497, 2014.

\bibitem{HalanayJMAA2015}
A.~Halanay and L.~Pandolfi.
\newblock Approximate controllability and lack of controllability to zero of
  the heat equation with memory.
\newblock {\em J. Math. Anal. Appl.}, 425(1):194--211, 2015.

\bibitem{HalanayDCDS-A2015}
A.~Halanay and L.~Pandolfi.
\newblock Noncontrollability for the {C}olemann-{G}urtin model in several
  dimensions.
\newblock {\em Discrete Contin. Dyn. Syst.}, pages 588--595, 2015.

\bibitem{HerraraGurtin1965}
I.~Herrera and M.~E. Gurtin.
\newblock A correspondence principle for viscoelastic wave propagation.
\newblock {\em Quart. Appl. Math.}, 22:360--364, 1965.

\bibitem{HopkinsonLONDON1876}
J.~Hopkinson.
\newblock The residual charge of the {L}eiden jar.
\newblock {\em Phil. Trans. Roy. Soc. Lond.}, 166:489--494, 1876.

\bibitem{HopkinsonLONDON1877}
J.~Hopkinson.
\newblock Residual charge of the {L}eiden jar---dielectric properties of
  different glasses.
\newblock {\em Phil. Trans. Roy. Soc. Lond.}, 167:569--626, 1877.

\bibitem{IvanovJMAA2009}
S.~Ivanov and L.~Pandolfi.
\newblock Heat equation with memory: Lack of controllability to rest.
\newblock {\em Journal of Mathematical Analysis and Applications},
  355(1):1--11, 2009.

\bibitem{JacobJMAA2007}
B.~Jacob and J.~R. Partington.
\newblock A resolvent test for admissibility of {V}olterra observation
  operators.
\newblock {\em J. Math. Anal. Appl.}, 332(1):346--355, 2007.

\bibitem{JacobBirgitJFA2010}
B.~Jacob, J.~R. Partington, and S.~Pott.
\newblock Weighted interpolation in {P}aley-{W}iener spaces and finite-time
  controllability.
\newblock {\em J. Funct. Anal.}, 259(9):2424--2436, 2010.

\bibitem{Jeffreys1924}
H.~Jeffreys.
\newblock {\em The Earth}.
\newblock {C}ambridge {U}niversity {P}ress, Cambridge, 1924.

\bibitem{JosephLIBRO1990}
D.~D. Joseph.
\newblock {\em Fluid dynamics of viscoelastic liquids}, volume~84.
\newblock Springer-Verlag, New York, 1990.

\bibitem{JosephPreziosi1989}
D.~D. Joseph and L.~Preziosi.
\newblock Heat waves.
\newblock {\em Rev. Modern Phys.}, 61(1):41--73, 1989.

\bibitem{JosephPreziosiAddendum1990}
D.~D. Joseph and L.~Preziosi.
\newblock Addendum to the paper: ``{H}eat waves'' [{R}ev.\ {M}odern {P}hys.\
  {\bf 61} (1989), no.\ 1, 41--73.
\newblock {\em Rev. Modern Phys.}, 62(2):375--391, 1990.

\bibitem{JungSIAMJCONTR2000}
M.~Jung.
\newblock Admissibility of control operators for solution families to
  {V}olterra integral equations.
\newblock {\em SIAM J. Control Optim.}, 38(5):1323--1333, 2000.

\bibitem{KabanikhinLorenziBOOK1999}
S.~I. Kabanikhin and A.~Lorenzi.
\newblock {\em Identification problems of wave phenomena: Theory and numerics}.
\newblock VSP, Utrecht, 1999.

\bibitem{KIMplate1992}
J.~U. Kim.
\newblock Control of a plate equation with large memory.
\newblock {\em Differential Integral Equations}, 5(2):261--279, 1992.

\bibitem{KimSIAM1993}
J.~U. Kim.
\newblock Control of a second-order integro-differential equation.
\newblock {\em SIAM J. Control Optim.}, 31(1):101--110, 1993.

\bibitem{KOHLRAUSCH1863}
F.~Kohlrausch.
\newblock Ueber dier elastiche narkwirkung bei dor torsion.
\newblock {\em Ann. Phys. Chem. (Pogg. Ann.)}, 119(3):337--368, 1992.

\bibitem{KolskyBOOK}
H.~Kolsky.
\newblock {\em Stress waves in solids}.
\newblock {D}over pubications, {I}nc., New York, 1963.

\bibitem{KrabsLeugeringSeidman1985}
W.~Krabs, G.~Leugering, and T.~I. Seidman.
\newblock On boundary controllability of a vibrating plate.
\newblock {\em Appl. Math. Optim.}, 13(3):205--229, 1985.

\bibitem{LagneseLionsLIBRO1988}
J.~Lagnese and J.-L. Lions.
\newblock {\em Modelling analysis and control of thin plates}.
\newblock Masson, Paris, 1988.

\bibitem{LasieckaKirchhoffCONF1989}
I.~Lasiecka.
\newblock Controllability of a viscoelastic {K}irchhoff plate.
\newblock In {\em Control and estimation of distributed parameter systems
  ({V}orau, 1988)}, volume~91 of {\em Internat. Ser. Numer. Math.}, pages
  237--247. Birkh\"auser, Basel, 1989.

\bibitem{LasieckaTriggianiCOSINE1981}
I.~Lasiecka and R.~Triggiani.
\newblock A cosine operator approach to modeling {$L_{2}(0,\,T;\ L_{2}(\Gamma
  ))$}---boundary input hyperbolic equations.
\newblock {\em Appl. Math. Optim.}, 7(1):35--93, 1981.

\bibitem{LavanyaBALACH}
R.~Lavanya and K.~Balachandran.
\newblock Controllability results of linear parabolic integrodifferential
  equations.
\newblock {\em Differential Integral Equations}, 21(9-10):801--819, 2008.

\bibitem{LeugeringFADING1984}
G.~Leugering.
\newblock Exact controllability in viscoelasticity of fading memory type.
\newblock {\em Applicable Anal.}, 18(3):221--243, 1984.

\bibitem{LeugeringONEDIMTHERMOVISCO1985}
G.~Leugering.
\newblock Boundary controllability in one-dimensional linear
  thermoviscoelasticity.
\newblock {\em J. Integral Equations}, 10(1-3, suppl.):157--173, 1985.

\bibitem{LeugeringBEAMAPPLANAL1986}
G.~Leugering.
\newblock Boundary controllability of a viscoelastic beam.
\newblock {\em Appl. Anal.}, 23(1-2):119--137, 1986.

\bibitem{LeugeringM2AS1986}
G.~Leugering.
\newblock Optimal controllability in viscoelasticity of rate type.
\newblock {\em Mathematical Methods in the Applied Sciences}, 8(1):368--386,
  1986.

\bibitem{LeugeringINTEGRODIFFERENTIAL1987}
G.~Leugering.
\newblock Exact boundary controllability of an integro-differential equation.
\newblock {\em Appl. Math. Optim.}, 15(3):223--250, 1987.

\bibitem{LeugeringLIQUIDM2AS1987}
G.~Leugering.
\newblock Time optimal boundary controllability of a simple linear viscoelastic
  liquid.
\newblock {\em Math. Methods Appl. Sci.}, 9(3):413--430, 1987.

\bibitem{LeugeringSTRING1989}
G.~Leugering.
\newblock Boundary controllability of a viscoelastic string.
\newblock In {\em Volterra integrodifferential equations in {B}anach spaces and
  applications ({T}rento, 1987)}, volume 190 of {\em Pitman Res. Notes Math.
  Ser.}, pages 258--270. Longman Sci. Tech., Harlow, 1989.

\bibitem{LeugeringDYNAstabSYSTEMS1989}
G.~Leugering.
\newblock On boundary feedback stabilization of a viscoelastic membrane.
\newblock {\em Dynam. Stability Systems}, 4(1):71--79, 1989.

\bibitem{LoretiSIAM2012}
P.~Loreti, L.~Pandolfi, and D.~Sforza.
\newblock Boundary controllability and observability of a viscoelastic string.
\newblock {\em SIAM J. Control Optim.}, 50(2):820--844, 2012.

\bibitem{LoretiCOMPRENDMATH2009}
P.~Loreti and D.~Sforza.
\newblock Exact reachability for second-order integro-differential equations.
\newblock {\em Comptes Rendus Mathematique}, 347(19-20):1153--1158, 2009.

\bibitem{LoretiJDE2010}
P.~Loreti and D.~Sforza.
\newblock Reachability problems for a class of integro-differential equations.
\newblock {\em Journal of Differential Equations}, 248(7):1711--1755, 2010.

\bibitem{Loreti2011261}
P.~Loreti and D.~Sforza.
\newblock Multidimensional controllability problems with memory.
\newblock {\em Operator Theory: Advances and Applications}, 216:261--274, 2011.

\bibitem{LoretiPARMA2016}
P.~Loreti and D.~Sforza.
\newblock Hidden regularity for wave equations with memory.
\newblock {\em Riv. Math. Univ. Parma (N.S.)}, 7(2):391--405, 2016.

\bibitem{LuJMATHPURAPPL2017}
Q.~L\"u, X.~Zhang, and E.~Zuazua.
\newblock Null controllability for wave equations with memory.
\newblock {\em Journal des Mathematiques Pures et Appliquees}, 108(4):500--531,
  2017.

\bibitem{MacCamyINTEGRODIFFERENTIAL1977}
R.~C. MacCamy.
\newblock An integro-differential equation with application in heat flow.
\newblock {\em Quart. Appl. Math.}, 35(1):1--19, 1977/78.

\bibitem{MacCamyMODEL1977}
R.~C. MacCamy.
\newblock A model for one-dimensional, nonlinear viscoelasticity.
\newblock {\em Quart. Appl. Math.}, 35(1):21--33, 1977/78.

\bibitem{MartinRosierSIAM2013}
P.~Martin, L.~Rosier, and P.~Rouchon.
\newblock Null controllability of the structurally damped wave equation with
  moving control.
\newblock {\em SIAM J. Control Optim.}, 51(1):660--684, 2013.

\bibitem{MaxwellPHILOSOPHICAL1867}
J.~C. Maxwell.
\newblock On the dynamical theory of gases.
\newblock {\em Phil. Trans. Roy. Soc. London}, 157:49--88, 1867.

\bibitem{Mitra2018575}
D.~Mitra, M.~Ramaswamy, and M.~Renardy.
\newblock Approximate controllability results for viscoelastic flows with
  infinitely many relaxation modes.
\newblock {\em Journal of Differential Equations}, 264(2):575--603, 2018.

\bibitem{MunozADVDE2003}
J.~E. Mu\~noz Rivera and M.~G. Naso.
\newblock Exact controllability in thermoelasticity with memory.
\newblock {\em Adv. Differential Equations}, 8(4):471--490, 2003.

\bibitem{MunozADVDE2004}
J.~E. Mu\~noz Rivera and M.~G. Naso.
\newblock Exact controllability for hyperbolic thermoelastic systems with large
  memory.
\newblock {\em Adv. Differential Equations}, 9(11-12):1369--1394, 2004.

\bibitem{NarainJosephLINEARIZED1983}
A.~Narain and D.~D. Joseph.
\newblock Linearized dynamics of shearing deformation perturbing rest in
  viscoelastic materials.
\newblock In {\em Equadiff 82 ({W}\"urzburg, 1982)}, volume 1017 of {\em
  Lecture Notes in Math.}, pages 476--507. Springer, Berlin, 1983.

\bibitem{NegrescuAPPLIED2016}
A.~Negrescu.
\newblock Controllability for the {N}eumann problem for the heat equation with
  memory.
\newblock {\em Appl. Comput. Math.}, 15(3):313--318, 2016.

\bibitem{PandolfiAMO2005}
L.~Pandolfi.
\newblock The controllability of the {G}urtin-{P}ipkin equation: A cosine
  operator approach.
\newblock {\em Applied Mathematics and Optimization}, 52(2):143--165, 2005.

\bibitem{PandolfiINTEQOPTH2009}
L.~Pandolfi.
\newblock Riesz systems and controllability of heat equations with memory.
\newblock {\em Integral Equations and Operator Theory}, 64(3):429--453, 2009.

\bibitem{PandolfiANNREVCONT2010}
L.~Pandolfi.
\newblock On-line input identification and application to active noise
  cancellation.
\newblock {\em Annual Reviews in Control}, 34(2):245--261, 2010.

\bibitem{PandolfiDCDSSerB2010}
L.~Pandolfi.
\newblock Riesz systems and moment method in the study of viscoelasticity in
  one space dimension.
\newblock {\em Discrete and Continuous Dynamical Systems - Series B},
  14(4):1487--1510, 2010.

\bibitem{PandolfiErratumAMO2011}
L.~Pandolfi.
\newblock Erratum to: {T}he controllability of the {G}urtin-{P}ipkin equation:
  a cosine operator approach.
\newblock {\em Appl. Math. Optim.}, 64(3):467--468, 2011.

\bibitem{PandolfiDCDSserS2011}
L.~Pandolfi.
\newblock Riesz systems, spectral controllability and a source identification
  problem for heat equations with memory.
\newblock {\em Discrete and Continuous Dynamical Systems - Series S},
  4(3):745--759, 2011.

\bibitem{PandolfiJMAA2013}
L.~Pandolfi.
\newblock Boundary controllability and source reconstruction in a viscoelastic
  string under external traction.
\newblock {\em J. Math. Anal. Appl.}, 407(2):464--479, 2013.

\bibitem{PandolfiEVOLEQCONTTH2013}
L.~Pandolfi.
\newblock Traction, deformation and velocity of deformation in a viscoelastic
  string.
\newblock {\em Evolution Equations and Control Theory}, 2(3):471--493, 2013.

\bibitem{PandolfiLIBRO2014}
L.~Pandolfi.
\newblock {\em Distributed systems with persistent memory. {C}ontrol and moment
  problems}.
\newblock SpringerBriefs in Control, Automation and Robotics,Springer, Cham,
  2014.

\bibitem{PandolfiJINTEQAPPL2015}
L.~Pandolfi.
\newblock Sharp control time for viscoelastic bodies.
\newblock {\em Journal of Integral Equations and Applications}, 27(1):103--136,
  2015.

\bibitem{PandolfiPARMA2016}
L.~Pandolfi.
\newblock Controllability for the heat equation with memory: a recent approach.
\newblock {\em Riv. Math. Univ. Parma (N.S.)}, 7(2):259--277, 2016.

\bibitem{PandolfiSIAM2017}
L.~Pandolfi.
\newblock Controllability of a viscoelastic plate using one boundary control in
  displacement or bending.
\newblock {\em SIAM J. Control Optim.}, 55(6):4092--4111, 2017.

\bibitem{PandolfiESAIM2017}
L.~Pandolfi.
\newblock Controllability of isotropic viscoelastic bodies of
  {M}axwell-{B}oltzmann type.
\newblock {\em ESAIM Control Optim. Calc. Var.}, 23(4):1649--1666, 2017.

\bibitem{PruessLIBRO1993}
J.~Pr\"uss.
\newblock {\em Evolutionary Integral Equations and Applications}.
\newblock Birkh\"auser, Basel, 1993.

\bibitem{RenardyREOLOGICAL1982}
M.~Renardy.
\newblock Some remarks on the propagation and nonpropagation of discontinuities
  in linearly viscoelastic liquids.
\newblock {\em Rheol. Acta}, 21(3):251--254, 1982.

\bibitem{RenardyCBMS2000}
M.~Renardy.
\newblock {\em Mathematical analysis of viscoelastic flows}, volume~73 of {\em
  CBMS-NSF Regional Conference Series in Applied Mathematics}.
\newblock Society for Industrial and Applied Mathematics (SIAM), Philadelphia,
  PA, 2000.

\bibitem{RenardySCL2005}
M.~Renardy.
\newblock Are viscoelastic flows under control or out of control?
\newblock {\em Systems and Control Letters}, 54(12):1183--1193, 2005.

\bibitem{Renardy2008}
M.~Renardy.
\newblock Chapter 5: {M}athematical {A}nalysis of {V}iscoelastic {F}luids.
\newblock {\em Handbook of Differential Equations: Evolutionary Equations},
  4:229--265, 2008.

\bibitem{RenardyHrusaLibro1987}
M.~Renardy, W.~J. Hrusa, and J.~A. Nohel.
\newblock {\em Mathematical problems in viscoelasticity}, volume~35 of {\em
  Pitman Monographs and Surveys in Pure and Applied Mathematics}.
\newblock Longman Scientific \& Technical, Harlow; John Wiley \& Sons, Inc.,
  New York, 1987.

\bibitem{RomanovShamev}
I.~Romanov and A.~Shamaev.
\newblock Exact controllability of the distributed system, governed by string
  equation with memory.
\newblock {\em J. Dyn. Control Syst.}, 19(4):611--623, 2013.

\bibitem{SaktivelJMAA2007}
K.~Sakthivel, K.~Balachandran, and R.~Lavanya.
\newblock Exact controllability of partial integrodifferential equations with
  mixed boundary conditions.
\newblock {\em J. Math. Anal. Appl.}, 325(2):1257--1279, 2007.

\bibitem{SaktivelCOMPUTERS2006}
K.~Sakthivel, K.~Balachandran, and S.~S. Sritharan.
\newblock Controllability and observability theory of certain parabolic
  integrodifferential equations.
\newblock {\em Comput. Math. Appl.}, 52(8-9):1299--1316, 2006.

\bibitem{ShubovSIAMcontr1997}
M.~A. Shubov, C.~F. Martin, J.~P. Dauer, and Boris~P. Belinskiy.
\newblock Exact controllability of the damped wave equation.
\newblock {\em SIAM J. Control Optim.}, 35(5):1773--1789, 1997.

\bibitem{SorianoCOMPLUTENSE1995}
J.~A. Soriano.
\newblock Exact controllability of the generalized telegraph equation.
\newblock {\em Rev. Mat. Univ. Complut. Madrid}, 8(2):459--493, 1995.

\bibitem{VlasovJOURNMATHSCI}
V.~V. Vlasov, Dzh. Vu, and G.~R. Kabirova.
\newblock Well-defined solvability and the spectral properties of abstract
  hyperbolic equations with aftereffect.
\newblock {\em J. Math. Sci. (N.Y.)}, 170(3):388--404, 2010.

\bibitem{VOLTERRArendiconti1909}
V~Volterra.
\newblock Sulle equazioni integro-differenziali.
\newblock {\em Rend. Acad. Naz. Lincei, Ser.~5}, XVIII:167--174, 1909.

\bibitem{VolterraACTAMATHEM1912}
V.~Volterra.
\newblock Sur les \'equations int\'egro-diff\'erentielles et leurs
  applications.
\newblock {\em Acta Math.}, 35:295--356, 1912.

\bibitem{VolterraLIBROfUNZINEA}
V.~Volterra.
\newblock {\em Le\c cons sur les fonctions de lignes}.
\newblock {G}authier-{V}illars,, Paris, 1913.

\bibitem{WangFeskanEVEQCONTRTH2017}
J.~Wang, M.~Feckan, and Y.~Zhou.
\newblock Approximate controllability of {S}obolev type fractional evolution
  systems with nonlocal conditions.
\newblock {\em Evol. Equ. Control Theory}, 6(3):471--486, 2017.

\bibitem{WangGuoFU2009}
J.-M. Wang, B.-Z. Guo, and M.-Y. Fu.
\newblock Dynamic behavior of a heat equation with memory.
\newblock {\em Math. Methods Appl. Sci.}, 32(10):1287--1310, 2009.

\bibitem{WesterlundTRANSDIEL1994}
S~Westerlund and L.~Ekstam.
\newblock Capacitors theory.
\newblock {\em IEEE Trans. Diel. Electr. Insul.}, 1(5):826--839, 1994.

\bibitem{YanCOMPLUTENSE1992}
J.~H. Yan.
\newblock Contr\^olabilit\'e exacte pour des syst\`emes \`a m\'emoire.
\newblock {\em Rev. Mat. Univ. Complut. Madrid}, 5(2-3):319--327, 1992.

\bibitem{FuZHANG-2}
J.~Yong and X.~Zhang.
\newblock Exact controllability of the heat equation with hyperbolic memory
  kernel.
\newblock In {\em Control theory of partial differential equations}, volume 242
  of {\em Lect. Notes Pure Appl. Math.}, pages 387--401. Chapman \& Hall/CRC,
  Boca Raton, FL, 2005.

\bibitem{ZuazuaTEMPSPETITparis1987}
E.~Zuazua.
\newblock Contr\^olabilit\'e exacte d'un mod\`ele de plaques vibrantes en un
  temps arbitrairement petit.
\newblock {\em C. R. Acad. Sci. Paris S\'er. I Math.}, 304(7):173--176, 1987.

\end{thebibliography}

 \enddocument